\documentclass[12pt]{article}
\usepackage{sw20elba}
\usepackage{amsmath,amssymb}

\newtheorem{theorem}{Theorem}[section]

\newtheorem{corollary}[theorem]{Corollary}

\newtheorem{definition}[theorem]{Definition}

\newtheorem{lemma}[theorem]{Lemma}

\newtheorem{proposition}[theorem]{Proposition}

\newenvironment{proof}[1][Proof]{\textbf{#1.} }{\ \rule{0.5em}{0.5em}}
\input{tcilatex}

\begin{document}

\author{Robert Erdahl and Konstantin Rybnikov}
\title{Voronoi-Dickson Hypothesis on Perfect Forms and L-types }
\date{November 19, 2001 }
\maketitle
\centerline{Short Version}
\begin{abstract}
George Voronoi (1908, 1909) introduced two important reduction methods for
positive quadratic forms: the reduction with perfect forms, and the
reduction with $L$-type domains, often called domains of Delaunay type. The
first method is important in studies of dense lattice packings of spheres.
The second method provides the key tools for finding the least dense lattice
coverings with equal spheres in lower dimensions. In his investigations
Voronoi heavily relied on that in dimensions less than 6 the partition of
the cone of positive quadratic forms into $L$-types refines the partition of
this cone into perfect domains. Voronoi conjectured implicitely and Dickson
(1972) explicitely that the $L$-partition is always a refinement of the
partition into perfect domains. This was proved for $n\leq 5$ (Voronoi,
Delaunay, Ryshkov, Baranovskii). We show that Voronoi-Dickson conjecture
fails already in dimension 6.
\end{abstract}


 \smallskip

\section{Introduction and main result\label{intro}}

Positive quadratic forms (referred to as PQFs) in $n$ indeterminate form a 
\emph{closed }cone $\frak{P}(n)$ of dimension $N=\frac{n(n+1)}{2}$ in $%
\mathbb{R}^{N}$, and this cone is the main object of study in our paper. The
interior of $\frak{P}(n)$ consists of positive definite forms of rank $n$.
We abbreviate positive \emph{definite} forms as PDQFs. PDQFs serve as
algebraic representations of \textit{point lattices}. There is a one-to-one
correspondence between isometry classes of $n$-lattices and integral
equivalence classes (i.e. with respect to $GL(n,\mathbb{Z})$-conjugation) of
PDQFs in $n$ indeterminates. For basic results of the theory of lattices and
PQFs and their applications see Ryshkov and Baranovskii (1978), Gruber and
Lekkerkerker (1987), Erd\"{o}s, Gruber and Hammer (1989), Conway and Sloane
(1999).

$GL(n,\mathbb{Z})$ acts pointwise on the space of quadratic forms $Sym(n,%
\mathbb{R)}\cong \mathbb{R}^{N}$. A \emph{polyhedral reduction partition} $R$
of $\frak{P}(n)$ is a partition of this cone into open \emph{convex
polyhedral} cones such that:

\begin{definition}
\begin{enumerate}
\item  it is invariant with respect to $GL(n,\mathbb{Z})$;

\item  there are finitely many inequivalent cones in this partition;

\item  for each cone $C$ of $R$ and any PQF $f$ in $n$ indeterminates, $f$
can be $GL(n,\mathbb{Z})$-equivalent to at most finitely many forms lying in 
$C$.
\end{enumerate}
\end{definition}

The \emph{partition into perfect cones} and the \emph{$L$-type partition}
(also referred to as the Voronoi partition of the 2nd kind, or the partition
into Voronoi reduction domains) are important polyhedral reduction
partitions of $\frak{P}(n)$. (Our usage of term \emph{domain }is lax; it
should be clear from the context whether we mean the whole arithmetic class,
or just \ one element of this class.) These partitions have been intensively
studied in geometry of numbers since times of Korkin, Zolotareff (1873) and
Voronoi (1908), and more recently in combinatorics (e.g. Deza et al.
(1997)), and algebraic geometry (e.g. Alexeev (1999a,b)). \emph{In most
previous works} (e.g. Voronoi (1908,1909), Ryshkov, Baranovskii (1976),
Dickson (1972)) \emph{the }$L$\emph{-type partition of }$P(n)$\emph{, or
sub-cones of }$P(n)$\emph{, was constructed by refining the perfect
partition.} It is not an exaggeration to say that in almost any systematic
study, except for Engels' computational investigations, $L$-types were
approached via perfect forms. For example, Voronoi started classifying
4-dimensional $L$-types by analyzing the Delaunay ($L$-)tilings of forms
lying in the 1st ($A_{n}$) and 2nd ($D_{n}$) perfect domains. The same route
was followed by Ryshkov and Baranovskii (1976). It was wideley believed that
the $L$-type parition is the refinement of the perfect parition, i.e. each
convex cone of the perfect partition is the union of finitely many conex
cones from the $L$-type partition. This conjecture is implicit in Voronoi's
memoirs (1908-1909), and explicit in Dickson (1972), where he showed that
the first perfect domain is the only perfect domain which coincided with an $%
L$-type domain.

In our paper we present the results of our study of the relashionship
between these two partitions for $n=6$. We have proved that for $n=6$ the $L$%
-type partition is NOT a refinement of the partition into perfect cones, as
was believed before. This note contains only a scketch of the proof of this
result, which will be later published in a longer journal paper (about
35-40pp).

\section{Perfect and $L$-type partitions$\label{Definitions}$}

\subsection{$L$-types}

\begin{definition}
Let $L$ be a lattice in $\mathbb{R}^{n}$. A convex polyhedron $P$ in $%
\mathbb{R}^{n}$ is called a Delaunay cell of $L$ with respect to a positive
quadratic form $f\mathbf{(x,x)}$ if:

\begin{enumerate}
\item  for each face $F$ of $P$ we have $conv(L\cap F)=F$;

\item  there is a quadric circumscribed about $P$, called the \emph{empty
ellipsoid} of $P$ when $f\mathbf{(x,x)}$ is positive definite, whose
quadratic form is $f\mathbf{(x,x)}$ (in case $rank\ f<n$, this quadric is an
elliptic cilinder);

\item  no points of $L$ lie inside the quadric circumscribed \ about $P$.
\end{enumerate}
\end{definition}

When $f=\sum_{i=1}^{n}x_{i}^{2}$, our definition coincides with the
classical definition of Delaunay cell in $\mathbb{E}^{n}$. Delaunay cells
form a convex face-to-face tiling of $L$ that is uniquely defined by $L$
(Delaunay, 1937). Two Delaunay cells are called \emph{homologous} is they
can be mapped to each other with a composition of a lattice translation and
a central inversion with respect to a lattice point.

\begin{definition}
PQFs $f_{1}$ and $f_{2}$ belong to the same convex $L$-domain if the
Delaunay tilings of $\mathbb{Z}^{n}$ with respect to $f_{1}$ and $f_{2}$ are
identical. $f_{1}$ and $f_{2}$ belong to the same $L$-type if these tilings
are equivalent with respect to $GL(n,\mathbb{Z)}$.
\end{definition}

The folowing proposition establishes the equivalence between the Delaunay's
definition of $L$-equivalence for lattices and the notion of $L$-equivalence
for arbitrary PQFs, which is introduced above.

\begin{proposition}
Positive definite forms $f$ and $g$ \ belong to the same $L$-type if the
corresponding lattices belong to the same\ $L$-type with respect to the form 
$\sum_{i=1}^{n}x_{i}^{2}$.
\end{proposition}

\begin{theorem}
(Voronoi) The parition of $\frak{P}(n)$ into $L$-types is a reduction
partition. Moreover, it is face-to-face.
\end{theorem}

The notions of Delauny tiling and $L$-type are important in the study of
extremal and group-theoretic properties of lattices. 
For example, the analysis of Delaunay cells in the Leech lattice conducted
by Conway, Sloane (e.g. \ see (1999)) and Borcherds showed that 23 ''deep
holes'' (Delaunay cells of radius equal to the covering radius of the
lattice) in the Leech lattice correspond to 23 even unimodular
24-dimensional lattices (Niemeier's list) that, in turn, give rise to 23
''gluing'' constructions of the Leech lattice from root lattices. Barnes and
Dickson (1967, 1968) and, later, in a geometric form, Delaunay et al. (1969,
1970) proved the following

\begin{theorem}
The \emph{closure} of any $N$-dimensional convex $L$-domain contains at most
one local minimum of the sphere covering density. The group of $GL(n,\mathbb{%
Z)}$-automorphisms of the domain maps this form to itself.
\end{theorem}

Using this approach, Delaunay, Ryshkov and Baranovskii (1963, 1976) found
the best lattice coverings in $\mathbb{E}^{4}$ and $\mathbb{E}^{5}.$ The
theory of $L$-types also has numerous connections to combinatorics and, in
particular, to cuts, hypermetrics, and regular graphs (see Deza et al.
(1997)). Recently, V. Alexeev (1999a,b) found exciting connections between
compactifications of moduli spaces of principally polirized abelian varities
and $L$-types.

\subsection{Perfect cones}

The $L$-type partition of $\frak{P}(n)$ is closely related to the theory of 
\emph{perfect forms} originated by Korkine and Zolotareff (1873). Let $f%
\mathbf{(x,x})$ be a PDQF. The \emph{arithmetic minimum } of $f\mathbf{(x,x}%
) $ is the minimum of this form on $\mathbb{Z}^{n}$. The integral vectors on
which this minimum is attained are called the representations of the
minimum, or the \emph{minimal vectors} of $f\mathbf{(x,x})$: these vectors
have the minimal length among all vectors of $\mathbb{Z}^{n}$ when $f\mathbf{%
(x,x})$ is used as the metrical form. Form $f\mathbf{(x,x})$ is called \emph{%
perfect} if it can be reconstructed up to scale from all representations of
its arithmetic minimum. In other words, a form $f\mathbf{(x,x})$ with the
arithmetic minimum $m$ and the set of minimal vectors $\{\mathbf{v}%
_{k}|\,k=1,...,2s\}$ is perfect if the system 
\begin{tabular}{|l|}
\hline
$\sum_{i,j=1}^{n}a_{ij}v_{k}^{i}v_{k}^{j}=m,$ \\ \hline
\end{tabular}
where $k=1,...,2s,$ has a unique solution $(a_{ij})$ in $Sym(n,\mathbb{R)}%
\cong \mathbb{R}^{N}$ (indeed, uniquness requires at least $n(n+1)$ minimal
vectors).

\begin{definition}
PQFs $f_{1}$ and $f_{2}$ belong to the same cone of the perfect partition if
they both can be written as \emph{strictly positive} linear combinations of
some subset of minimal vectors of a perfect form $\phi $. $f_{1}$ and $f_{2}$
belong to the same perfect type if there is $f_{1}^{\prime }$, equivalent to 
$f_{1}$, such that $f_{1}^{\prime }$ and $f_{2}$ belong to the same cone of
the perfect partition.
\end{definition}

\begin{theorem}
(Voronoi) The partition of $\frak{P}(n)$ into perfect domains is a reduction
partition. Moreover, it is face-to-face. Each 1-dimensional cone of this
partition lies on $\partial \frak{P}(n)$.
\end{theorem}

Perfect forms play an important role in lattice sphere packings. Voronoi's
theorem (1908) says that if a form is\emph{\ extreme}---i.e., a maximum of
the packing density---it must also be perfect (see Coxeter (1951), Conway,
Sloane (1988) for the proof). The notion of \emph{eutactic form} arises in
the study of the dense lattice sphere packings and is directly related to
the notion of perfect form. The reciprocal\ of $f(\mathbf{x,x})$ is a form
whose Gramm matrix is the inverse of the Gramm matrix of $f(\mathbf{x,x}).$
The dual form is normally denoted by $f^{\ast }(\mathbf{x,x}).$ A form $f(%
\mathbf{x,x})$ is called \textit{eutactic }if the dual form $f^{\ast }(%
\mathbf{x,x})$ can be written as $\dsum_{k=1}^{s}\alpha _{k}(\mathbf{v}%
_{k}\bullet \mathbf{x)}^{2},$ where $\{\mathbf{v}_{k}|\,k=1,...,s\}$ is the
set of mutually non-collinear minimal vectors of $f(\mathbf{x,x}),$ and $%
\alpha _{k}>0$.

\begin{theorem}
(Voronoi) A\ form $f(\mathbf{x,x})$ is a maximum of the sphere packing
density if and only if $f(\mathbf{x,x})$ is perfect and eutactic.
\end{theorem}

Voronoi gave an algorithm finding all perfect domains for given $n$. This
algorithm is known as Voronoi's reduction with perfect forms. \ For the
computational analysis of his algorithm and its improvements see Martinet
(1996). The perfect forms and the incidence graphs of perfect paritions of $%
\frak{P}(n)$ have have been completely described for $n\leq 7$

\subsection{The relashionship in low dimensions. The case of $n=6$.}

Voronoi (1908-09) proved that for $n=2,3$ the $L$-partition and the perfect
partition of $\frak{P}(n)$ coincide. The perfect facet $D_{4}$ (the 2nd
perfect form in 4 variables) exemplifies a new pattern in the relation of
these partitions. Namely, the facet $D_{4}$ is decomposed into a number of
simplicial $L$-type domains like a pie:\ this decomposition consists of the
cones with apex at the affine center of this facet over the $(N-2)$-faces.
These simplexes are $L$-type domains of two arithmetic types: type I is
adjacent to the the perfect/$L$-type domain of $A_{4},$ type II is adacent
to an arithmetically equivalent $L$-type domain (also type II, indeed) from
the $L$-subdivision of the adjacent $D_{4}$ domain (for details see Delaunay
et al. (1963, 1968)).

Voronoi also proved that for $n=4$ the tiling of $\frak{P}(n)$ with $L$-type
domains refines the partition of this cone into perfect domains. Ryshkov and
Baranovskii (1975) proved the refinement hypothesis for $n=5$. In his paper
of 1972 T.J. Dickson proved that the perfect domain of $%
A_{n}$, also called the \textit{first perfect form} after Korkine and
Zolotareff (1873), is the only perfect domain that is also an $L$-type
domain; he was also first to explicitly mention the common believe in 
\emph{Voronoi's refinement hypothesis}.

\section{Metrical forms for the lattices $E_{6\text{ }}$and $E_{6}^{\ast }%
\label{FORMS}$}

Consider the following symetric sets of vectors in $\mathbb{Z}^{6}$. 
\begin{gather*}
\mathcal{P}_{1}=\{\pm \lbrack -3;2^{5}],\pm \lbrack 2;-2;-1^{4}],\pm \lbrack
1;0;-1^{4}]\} \\
\mathcal{P}_{2}=\{\pm \lbrack 2;-1;-2,-1^{3}]\times 4,\pm \lbrack
1;-1;0,-1^{3}]\times 4,\pm \lbrack 1;0^{2},-1^{3}]\times 10,\pm \lbrack
0;1,0^{4}]\times 5,\pm \lbrack 2;-1^{5}]\} \\
\mathcal{P}_{3}=\{\pm \lbrack 0;0;1,-1,0^{2}]\times 6,\pm \lbrack
1;0;-1^{2},0^{2}]\times 6\}
\end{gather*}
\ Here we use a short-hand notation for families of vectors obtained from
some $n$-vector by all circular permutations of selected subsets of its
components: (1) $m^{k}$ stands for $k$ consecutive $m$'s, (2) square
brackets $[a_{1}...a_{n}]$ are used to denote all vectors that can be
obtained from vector $(a_{1}...a_{n})$ by circular permutations in strings
of symbols that are separated by commas and bordered on the sides by
semicolons and/or brackets, (3) numbers between semicolumns and/or brackets
are not permuted.

The two sets of \emph{perfect vectors} $\mathcal{P}_{E_{6}^{\ast }}=\mathcal{%
P}_{1}\cup \mathcal{P}_{2}$ and $\mathcal{P}_{E_{6}}=\mathcal{P}_{2}\cup 
\mathcal{P}_{3}$ are the minimal vectors for the metrical forms $\pi
_{E_{6}^{\ast }}$ and $\pi _{E_{6}}$, with \emph{arithmetic minimum} $m$ and
coefficient matrices 
\begin{equation*}
\mathbf{P}_{E_{6}}=\frac{m}{2}\left[ 
\begin{array}{rrrrrr}
\;8 & \;1 & \;3 & \;3 & \;3 & \;3 \\ 
1 & 2 & 0 & 0 & 0 & 0 \\ 
3 & 0 & 2 & 1 & 1 & 1 \\ 
3 & 0 & 1 & 2 & 1 & 1 \\ 
3 & 0 & 1 & 1 & 2 & 1 \\ 
3 & 0 & 1 & 1 & 1 & 2
\end{array}
\right] ,\qquad \mathbf{P}_{E_{6}^{\ast }}=\frac{m}{4}\left[ 
\begin{array}{rrrrrr}
16 & \;5 & \;5 & \;5 & \;5 & \;5 \\ 
5 & 4 & 1 & 1 & 1 & 1 \\ 
5 & 1 & 4 & 1 & 1 & 1 \\ 
5 & 1 & 1 & 4 & 1 & 1 \\ 
5 & 1 & 1 & 1 & 4 & 1 \\ 
5 & 1 & 1 & 1 & 1 & 4
\end{array}
\right] .
\end{equation*}
The metrical forms $\pi _{E_{6}^{\ast }}$ and $\pi _{E_{6}}$\ are integrally
equivalent to forms $\phi _{4}$ and $\phi _{2}$ from Barnes's table (1957).
That $\pi _{E_{6}^{\ast }}$ and $\pi _{E_{6}}$ are forms for the root
lattices $E_{6}^{\ast }$ and $E_{6}^{\ast }$ is established by the numbers
of perfect vectors, which are given by $|\mathcal{P}_{E_{6}^{\ast }}|=54$
and $|\mathcal{P}_{E_{6}}|=72$. The. As shown below, if this scale parameter
is set equal to $\sqrt{8/3}$, then, the geometric lattices corresponding to $%
\pi _{E_{6}}$, $\pi _{E_{6}^{\ast }}$ are dual lattices.

\paragraph{Perfect vectors and perfect domains.}

Let $V\subset \mathbb{Z}^{d}$ be centrally symmetric. \ Then define the 
\emph{domain of} $V$ to be given by 
\begin{equation*}
\Phi (V)=\{\varphi (\mathbf{x})=\sum_{\mathbf{v\in }V^{+}}\omega _{\mathbf{v}%
}(\mathbf{v}\cdot \mathbf{x})^{2}|\omega _{\mathbf{v}}\geq 0\};
\end{equation*}
as indicated by the superscript $+$, this summation is over an \emph{%
oriented subset} of vectors -- one vector from each pair of opposites in the
symmetric set $V$. \ The domains $\Phi (\mathcal{P}_{E_{6}}),\Phi (\mathcal{P%
}_{E_{6}^{\ast }})$ are \emph{perfect cones} of types$\ E_{6}$ and $%
E_{6}^{\ast }$, respectively; they are defined on the perfect vectors for
the perfect forms $\pi _{E_{6}}$ and $\pi _{E_{6}^{\ast }}$.

\begin{proposition}
$\Phi (\mathcal{P}_{E_{6}}\cap \mathcal{P}_{E_{6}^{\ast }})$ is a facet of
both$\ \Phi (\mathcal{P}_{E_{6}})$ and $\Phi (\mathcal{P}_{E_{6}^{\ast }})$.
\end{proposition}

\begin{proof}
If forms $\pi $, $\varphi $ have coefficient matrices $\mathbf{P}$, $\mathbf{%
F}$, define the scalar product $\left\langle \pi ,\varphi \right\rangle
:=trace(\mathbf{PF)}$. \ Consider the form $\pi _{E_{6}^{\ast }}-\pi _{E_{6}}
$ and the rank one form $\varphi _{\mathbf{p}}(\mathbf{x})=(\mathbf{p\cdot x)%
}^{2}$. \ Then, $\left\langle \pi _{E_{6}^{\ast }}-\pi _{E_{6}},\varphi _{%
\mathbf{p}}\right\rangle =trace(\mathbf{P}_{E_{6}^{\ast }}-\mathbf{P}%
_{E_{6}})\mathbf{pp}^{T}=\mathbf{p}^{T}(\mathbf{P}_{E_{6}^{\ast }}-\mathbf{P}%
_{E_{6}})\mathbf{p}=\pi _{E_{6}^{\ast }}(\mathbf{p})-\pi _{E_{6}}(\mathbf{p})
$. \ If $\mathbf{p}\in \mathcal{P}_{E_{6}^{\ast }}$, then $\pi _{E_{6}^{\ast
}}(\mathbf{p})-\pi _{E_{6}}(\mathbf{p})=m-\pi _{E_{6}}(\mathbf{p})\leq 0$,
with equality if and only if $\mathbf{p}\in \mathcal{P}_{E_{6}}\cap \mathcal{%
P}_{E_{6}^{\ast }}$; if $\mathbf{p}\in \mathcal{P}_{E_{6}}$, then $\pi
_{E_{6}^{\ast }}(\mathbf{p})-\pi _{E_{6}}(\mathbf{p})=\pi _{E_{6}^{\ast }}(%
\mathbf{p})-m\geq 0$, with equality if and only if $\mathbf{p}\in \mathcal{P}%
_{E_{6}}\cap \mathcal{P}_{E_{6}^{\ast }}$. \ It follows that the hyperplane
with equation $\left\langle \pi _{E_{6}^{\ast }}-\pi _{E_{6}},\varphi
\right\rangle =0$ separates $\Phi (\mathcal{P}_{E_{6}})$ and $\Phi (\mathcal{%
P}_{E_{6}^{\ast }})$, and that $\Phi (\mathcal{P}_{E_{6}}\cap \mathcal{P}%
_{E_{6}^{\ast }})$ is a face of both these perfect domains. \ \ 

That $\Phi (\mathcal{P}_{E_{6}}\cap \mathcal{P}_{E_{6}^{\ast }})$ is a facet
follows by showing that the linear span of the forms $\varphi _{\mathbf{p}}(%
\mathbf{x})=(\mathbf{p\cdot x)}^{2}$, $\mathbf{p}\in \mathcal{P}_{E_{6}}\cap 
\mathcal{P}_{E_{6}^{\ast }}$ has co-dimension one in the linear space of
metrical forms. (We omit this argument.)
\end{proof}

\paragraph{Eutactic forms.}

The forms 
\begin{equation*}
\varphi _{E_{6}}(\mathbf{x})=\frac{m}{12}\sum_{\mathbf{p\in }\mathcal{P}%
_{E_{6}^{\ast }}^{+}}(\mathbf{p\cdot x})^{2},\qquad \varphi _{E_{6}^{\ast }}(%
\mathbf{x})=\frac{m}{16}\sum_{\mathbf{p\in }\mathcal{P}_{E_{6}}^{+}}(\mathbf{%
p\cdot x})^{2},
\end{equation*}
lie on the \emph{central rays} of the perfect cones $\mathcal{\Phi }%
_{E_{6}^{\ast }}$, $\mathcal{\Phi }_{E_{6}}$, and have coefficient matrices
given by 
\begin{equation*}
\mathbf{F}_{E_{6}}=\frac{m}{2}\left[ 
\begin{array}{rrrrrr}
8 & -5 & -5 & -5 & -5 & -5 \\ 
-5 & 4 & 3 & 3 & 3 & 3 \\ 
-5 & 3 & 4 & 3 & 3 & 3 \\ 
-5 & 3 & 3 & 4 & 3 & 3 \\ 
-5 & 3 & 3 & 3 & 4 & 3 \\ 
-5 & 3 & 3 & 3 & 3 & 4
\end{array}
\right] ,\text{ }\mathbf{F}_{E_{6}^{\ast }}=\frac{m}{4}\left[ 
\begin{array}{rrrrrr}
10 & -5 & -6 & -6 & -6 & -6 \\ 
-5 & 4 & 3 & 3 & 3 & 3 \\ 
-6 & 3 & 6 & 3 & 3 & 3 \\ 
-6 & 3 & 3 & 6 & 3 & 3 \\ 
-6 & 3 & 3 & 3 & 6 & 3 \\ 
-6 & 3 & 3 & 3 & 3 & 6
\end{array}
\right] .
\end{equation*}
These forms are related to the original $\pi _{E_{6}}$, $\pi _{E_{6}^{\ast
}} $ by the formulas $\varphi _{E_{6}}(\mathbf{x})=\pi _{E_{6}}(\mathbf{U}%
_{\ast }\mathbf{x})$, $\varphi _{E_{6}^{\ast }}(\mathbf{x})=\pi
_{E_{6}^{\ast }}(\mathbf{U}_{\ast }\mathbf{x})$, where $\mathbf{U}_{\ast
}\in GL(6,\mathbb{Z})$ is given by 
\begin{equation*}
\mathbf{U}_{\ast }\mathbf{=}\left[ 
\begin{array}{rrrrrr}
2 & -2 & 0 & -1 & -1 & -1 \\ 
-2 & 2 & 1 & 1 & 1 & 1 \\ 
0 & 1 & -1 & 0 & 0 & 0 \\ 
-1 & 1 & 0 & 1 & 0 & 0 \\ 
-1 & 1 & 0 & 0 & 1 & 0 \\ 
-1 & 1 & 0 & 0 & 0 & 1
\end{array}
\right] .\qquad \mathbf{U}_{\ast }^{-1}\mathbf{=}\left[ 
\begin{array}{rrrrrr}
0 & 1 & 1 & -1 & -1 & -1 \\ 
1 & 1 & 1 & 0 & 0 & 0 \\ 
1 & 1 & 0 & 0 & 0 & 0 \\ 
-1 & 0 & 0 & 0 & -1 & -1 \\ 
-1 & 0 & 0 & -1 & 0 & -1 \\ 
-1 & 0 & 0 & -1 & -1 & 0
\end{array}
\right]
\end{equation*}
Being arithmetically equivalent, $\varphi _{E_{6}}$, $\pi _{E_{6}}$ are
alternate metrical forms for the same geometric lattice, and similarly for $%
\varphi _{E_{6}^{\ast }}$, $\pi _{E_{6}^{\ast }}$. We will refer to the
minimal vectors of these forms as \emph{short}. \ The short vectors for $%
\varphi _{E_{6}}$ and $\varphi _{E_{6}^{\ast }}$\ are related to the perfect
vectors for $\pi _{E_{6}}$ and $\pi _{E_{6}^{\ast }}$\ by the formulas: $%
\mathcal{S}_{E_{6}}=\mathcal{S}_{2}\cup \mathcal{S}_{3}$, $\mathcal{S}%
_{E_{6}^{\ast }}=\mathcal{S}_{1}\cup \mathcal{S}_{2}$, where $\mathcal{S}%
_{1}=\mathbf{U}_{\ast }^{-1}(\mathcal{P}_{1})$, $\mathcal{S}_{2}=\mathbf{U}%
_{\ast }^{-1}(\mathcal{P}_{2})$, $\mathcal{S}_{3}=\mathbf{U}_{\ast }^{-1}(%
\mathcal{P}_{3})$. \ Explicit form for the sets $\mathcal{S}_{1}$, $\mathcal{%
S}_{2}$, $\mathcal{S}_{3}$ can be obtained by a direct calculation. 
\begin{gather*}
\mathcal{S}_{1}=\{\pm \lbrack 2,-1,1^{4}],\pm \lbrack -2,0,-1^{4}],\pm
\lbrack 0,1,0^{4}]\} \\
\mathcal{S}_{2}=\{\pm \lbrack 0,1;-1,0^{3}]\times 4,\pm \lbrack
2,0;0,1^{3}]\times 4,\pm \lbrack 1;1^{2},0^{3}]\times 10,\pm \lbrack
1;1,0^{4}]\times 5,\pm \lbrack 3;1^{5}]\} \\
\mathcal{S}_{3}=\{\pm \lbrack 0,0;1,-1,0^{2}]\times 6,\pm \lbrack
2,1;1^{2},0^{2}]\times 6\}
\end{gather*}

The four coefficient matrices satisfy the relations $\mathbf{F}_{E_{6}}%
\mathbf{P}_{E_{6}^{\ast }}=\mathbf{F}_{E_{6}^{\ast }}\mathbf{P}_{E_{6}}=%
\frac{3}{8}m^{2}\mathbf{I}$, so that when $m=\sqrt{8/3}$, $\mathbf{F}%
_{E_{6}}=(\mathbf{P}_{E_{6}^{\ast }})^{-1}$, $\mathbf{F}_{E_{6}^{\ast }}=(%
\mathbf{P}_{E_{6}})^{-1}$. Under these circumstances the pair of forms $%
\varphi _{E_{6}}$, $\pi _{E_{6}^{\ast }}$, and the pair $\varphi
_{E_{6}^{\ast }}$, $\pi _{E_{6}}$, are in duality: $\varphi _{E_{6}}=\pi
_{E_{6}^{\ast }}^{\circ }$, $\pi _{E_{6}^{\ast }}=\varphi _{E_{6}}^{\circ }$%
, $\varphi _{E_{6}^{\ast }}=\pi _{E_{6}}^{\circ }$, $\pi _{E_{6}}=\varphi
_{E_{6}^{\ast }}^{\circ }$. \ Dual forms correspond to dual lattices, so the
single geometric lattice corresponding to the forms $\varphi _{E_{6}}$, $\pi
_{E_{6}}$ is dual to the geometric lattice corresponding to $\varphi
_{E_{6}^{\ast }}$, $\pi _{E_{6}^{\ast }}$. This pair of geometric lattices
is the pair of root lattices $E_{6}$ and $E_{6}^{\ast }$.

The above representations for $\varphi _{E_{6}}$, $\varphi _{E_{6}^{\ast }}$
can equally well be considered as representations for the dual forms $\pi
_{E_{6}^{\ast }}^{\circ }=\varphi _{E_{6}}$, $\pi _{E_{6}}^{\circ }=\varphi
_{E_{6}^{\ast }}$, which shows that $\pi _{E_{6}^{\ast }}^{\circ }$, $\pi
_{E_{6}}^{\circ }$ lie interior to the domains $\mathcal{\Phi }_{E_{6}^{\ast
}}$, $\mathcal{\Phi }_{E_{6}}$ determined, respectively, by $\pi
_{E_{6}^{\ast }}$, $\pi _{E_{6}}$. Since $\pi _{E_{6}}=\varphi _{E_{6}^{\ast
}}^{\circ }$, $\pi _{E_{6}^{\ast }}=\varphi _{E_{6}}^{\circ }$, the original
forms have \emph{eutaxy }(see Section 2)\emph{\ }representations identicle
to those for $\varphi _{E_{6}}$, $\varphi _{E_{6}^{\ast }}$, 
\begin{equation*}
\pi _{E_{6}}(\mathbf{x})=\frac{m}{12}\sum_{\mathbf{s\in }\mathcal{S}%
_{E_{6}^{\ast }}^{+}}(\mathbf{s}\cdot \mathbf{x})^{2},\text{\qquad }\pi
_{E_{6}^{\ast }}(\mathbf{x})=\frac{m}{16}\sum_{\mathbf{s\in }\mathcal{S}%
_{E_{6}}^{+}}(\mathbf{s}\cdot \mathbf{x})^{2};
\end{equation*}
the summation again is over an oriented subset of minimal vectors for the
reciprocal forms, respectively, the forms $\varphi _{E_{6}^{\ast }}$, $%
\varphi _{E_{6}}$.

The symmetrical arrangement provided by the domains, $\Phi _{E_{6}^{\ast }}$%
, $\Phi _{E_{6}}$, is remarkable -- it is a rare occurance that forms on the
central rays of adjacent domains correspond to dual geometric lattices. \
While the lattices are dual, the forms $\varphi _{E_{6}}$, $\varphi
_{E_{6}^{\ast }}$ are not. \ These forms fit into the following intricate
pattern: $\varphi _{E_{6}}$ is dual to $\pi _{E_{6}^{\ast }}$, which in turn
is arithmetically equivalent to $\varphi _{E_{6}^{\ast }}$ and, $\varphi
_{E_{6}^{\ast }}$ is dual to $\pi _{E_{6}}$, which in turn is arithmetically
equivalent to $\varphi _{E_{6}}$. This pattern \ is known since the times of
Korkine and Zolotareff (1873). See also Coxeter (1951). Barnes (1957) showed
that there is only one arithmetic type of wall between $E_{6}$ and $%
E_{6}^{\ast }$.

\subsubsection{ The invariance groups for the forms $\protect\varphi %
_{E_{6}} $ and \ $\protect\varphi _{E_{6}^{\ast }}$.}

The point group for the root lattice $E_{6}$ is the product of the two
element group generated by central inversion and the reflection group $E_{6}$%
, and has order $2^{8}3^{4}5$. \ There is a representation $\mathcal{G}%
_{E_{6}}\subset $ $GL(6,\mathbb{Z)}$ that is the invariance group of the
perfect form $\varphi _{E_{6}}$, with the following characterization: $%
\mathbf{g\in }\mathcal{G}_{E_{6}}$ if and only if $\mathcal{\varphi }%
_{E_{6}}(\mathbf{gx)=}\mathcal{\varphi }_{E_{6}}(\mathbf{x)}$. \ There is a
second representation $\mathcal{G}_{E_{6}^{\ast }}$ that is the invariance
group of the form $\varphi _{E_{6}^{\ast }}$. \ 

The dual group (representation) $\mathcal{G}_{E_{6}}^{\circ }$ is the
invariance group of the perfect form $\pi _{E_{6}^{\ast }}$, and the dual
group (representation) $\mathcal{G}_{E_{6}^{\ast }}^{\circ }$ is the
invariance group of the form $\pi _{E_{6}}$. \ The subgroup $\mathcal{G}%
_{E_{6}}^{\circ }$ is defined by 
\begin{equation*}
\mathcal{G}_{E_{6}}^{\circ }=\{(\mathbf{g}^{-1})^{T}|\mathbf{g\in }\mathcal{G%
}_{E_{6}}\},
\end{equation*}
\ where $\mathbf{g}^{\circ }\mathbf{=}(\mathbf{g}^{-1})^{T}$ is the matrix
for the dual transformation; there is a similar definition for the subgroup $%
\mathcal{G}_{E_{6}^{\ast }}^{\circ }$. \ From the arithmetic equivalences $%
\mathcal{P}_{E_{6}^{\ast }}=\mathbf{U}_{\ast }\mathcal{S}_{E_{6}^{\ast }}$
and $\mathcal{P}_{E_{6}}=\mathbf{U}_{\ast }\mathcal{S}_{E_{6}}$, it follows
that $\mathcal{G}_{E_{6}}^{\circ }=\mathbf{U}_{\ast }\mathcal{G}%
_{E_{6}^{\ast }}\mathbf{U}_{\ast }^{-1}$ and $\mathcal{G}_{E_{6}^{\ast
}}^{\circ }=\mathbf{U}_{\ast }\mathcal{G}_{E_{6}}\mathbf{U}_{\ast }^{-1}$,
where $\mathbf{U}_{\ast }\in GL(6,\mathbb{Z})$ was described above. \ 

\begin{proposition}
\label{Invariance.P}The subgroups $\mathcal{G}_{E_{6}},\mathcal{G}%
_{E_{6}^{\ast }},\mathcal{G}_{E_{6}}^{\circ },\mathcal{G}_{E_{6}^{\ast
}}^{\circ }\subset GL(6,\mathbb{Z)}$ are the full invariance groups of the
respective sets of vectors $\mathcal{S}_{E_{6}},\mathcal{S}_{E_{6}^{\ast }},%
\mathcal{P}_{E_{6}^{\ast }},\mathcal{P}_{E_{6}}$, and act transitively on
these sets.
\end{proposition}

The remaining sections will be devoted to the proof of the \emph{main theorem%
} that provides for a counterexample to the refinement hypothesis:

\begin{theorem}
(main)\label{Main}The segment $t\varphi _{E_{6}}+(1-t)\varphi _{E_{6}^{\ast
}},\;0\leq t\leq 1$ has forms of 5 $L$-types: $t=0,\;0<t<1/2,\;t=1/2,%
\;1/2<t<1,\;t=1$. The wall between the perfect domains $\Phi (\mathcal{P}%
_{E_{6}})$ and $\Phi (\mathcal{P}_{E_{6}^{\ast }})$ crosses this segment at
point $\frac{2}{5}\varphi _{E_{6}}+\frac{3}{5}\varphi _{E_{6}^{\ast }}$.
\end{theorem}

\section{$L$-types}

\subsection{Commensurate Delaunay Tilings}

A convex polyhedron $P$ in $\mathbb{R}^{n}$ is called a $\mathbb{Z}^{n}$%
-polyhedron if any face $F$\ of $P$ is $conv(aff\;F$ $\cap $ $\mathbb{Z}%
^{n}).$ A $\mathbb{Z}^{n}$-tiling is a face-to-face tiling of $\mathbb{R}%
^{n} $ by convex $\mathbb{Z}^{n}$-polyhedra. If $T_{1},$ $T_{2}$ are two
convex tilings of $\mathbb{R}^{n}$, one can define the intersection tiling
of $T_{1} $ and $T_{2}:$ the (open) tiles of this tiling are all non-empty
intersections of (open) tiles of $T_{1}$ and $T_{2}.$ The intersection
tiling of two $\mathbb{Z}^{n}$-tilings is not always a $\mathbb{Z}^{n}$%
-tiling; it is a $\mathbb{Z}^{n}$-tiling if and only if the vertex set of
the intersection tiling is $\mathbb{Z}^{n}$.

\begin{definition}
Let $P_{1}$ and $P_{2}$ be two $\mathbb{Z}^{n}$-polyhedra. They are called
commensurate if all the faces of $P_{1}\cap P_{2}$ are $\mathbb{Z}^{n}$%
-polyhedra.
\end{definition}

In particular, two $\mathbb{Z}^{n}$-\emph{polytopes} $P_{1}$ and $P_{2}$ are
commensurate if the vertex set of $P_{1}\cap P_{2}$ belongs to $\mathbb{Z}%
^{n}.$

\begin{definition}
Let $T_{1},$ $T_{2}$ be two $\mathbb{Z}^{n}$-tilings. They are called
commensurate\ if their intersection tiling consist of $\mathbb{Z}^{n}$%
-polyhedra only.
\end{definition}

In particular, the vertex set of the intersection tiling of $T_{1}$ and $%
T_{2}$ must be a subset of $\mathbb{Z}^{n}$ (possibly empty).

\begin{lemma}
\label{Test on Commensurability of 2 cells}Assume that lattice tilings $%
\mathcal{D}(\varphi _{1}),\mathcal{D}(\varphi _{2})$ of $\mathbb{Z}^{n}$ are
Delaunay with respect to the PQFs $\varphi _{1}$ and $\varphi _{2}$. \ Also
assume that $C_{1}$, $C_{2}$ are closed Delaunay cells for $\mathcal{D}%
(\varphi _{1})$ and $\mathcal{D}(\varphi _{2})$. \ If $C_{1}\cap C_{2}\cap 
\mathbb{Z}^{d}\neq \emptyset $, then $C=conv(C_{1}\cap C_{2}\cap \mathbb{Z}%
^{d})$ is a closed Delaunay cell for all intermediate forms $\varphi
_{t}=(1-t)\varphi _{1}+t\varphi _{2}$, where $0<t<1$.
\end{lemma}

\begin{proof}
\label{ComTest} Since $C_{1}$ is Delaunay, there is a scalars $c_{1}$ and
vector $\mathbf{p}_{1}$ so that $f_{1}(\mathbf{x})=c_{1}+\mathbf{p}_{1}\cdot 
\mathbf{x+}\varphi _{1}(\mathbf{x})$ is zero-valued on the vertices of $%
C_{1} $, but has positive values on all other elements of $\mathbb{Z}^{d}$.
\ Therefore, the equation $f_{1}(\mathbf{x})=0$ is for an empty ellipsoid
circumscribing $C_{1}$. \ There is a similar function for $C_{2}$, so that $%
f_{2}(\mathbf{x})=c_{2}+\mathbf{p}_{2}\cdot \mathbf{x+}\varphi _{2}(\mathbf{x%
})=0$ is the equation of an empty ellipsoid circumscribing $C_{2}$. \ For $%
0<t<1$, the function $f_{t}(\mathbf{x})=(1-t)f_{1}(\mathbf{x})+tf_{2}(%
\mathbf{x})$ is zero-valued on the vertices of $C$, and positive on all
other elements of $\mathbb{Z}^{d}$. \ Therefore, $C$ is circumscribed by the
empty ellipsoid with equation $f_{t}(\mathbf{x})=0$, and is Delaunay with
respect to the form $\varphi _{t}=(1-t)\varphi _{1}+t\varphi _{2}$.
\end{proof}

\begin{lemma}
\label{Commense} If Delaunay tilings $\mathcal{D}(f)$ and $\mathcal{D}(g)$
of $\mathbb{Z}^{n}$ with respect to PQFs $f$ and $g$ are commensurate, then
their intersection $\mathbb{Z}^{n}$-tiling is Delaunay with respect to $%
\alpha f+\beta g$ for any $\alpha ,\beta >0.$ Conversly, if for all $0<t<1$
the Delaunay tiling of $\mathbb{Z}^{n}$ \ with respect to $tf+(1-t)g$ is the
same, then Delaunay tilings $\mathcal{D}(f)$ and $\mathcal{D}(g)$ of $%
\mathbb{Z}^{n}$ are commensurate.
\end{lemma}

\begin{proof}
$\Longrightarrow )$ Since $\mathcal{D}(f)$ and $\mathcal{D}(g)$ are
commensurate, any cell $C$ of the intersection tiling of $\mathcal{D}(f)$
and $\mathcal{D}(g)$ can be written as $C=conv(C_{1}\cap C_{2}\cap \mathbb{Z}%
^{d}),$ where $C_{1}$, $C_{2}$ are closed Delaunay cells for $\mathcal{D}%
(\varphi _{1})$ and $\mathcal{D}(\varphi _{2})$. By the above lemma, $C$ is
Delaunay for $\alpha f+\beta g$ for any $\alpha ,\beta >0.$

$\Longleftarrow )$ If $C$ is a Delaunay cell with respect to $tf+(1-t)g$,
then, by a standard continuity argument, there are $\mathbb{Z}^{n}$%
-polyhedra $C_{f}$,$C_{g}$ which are Delaunay relative to $f$ and $g$
respectively, such that $C\subseteq C_{f}$ and $C\subseteq C_{g}.$
Therefore, $C=C_{f}\cap C_{g}=conv(C_{f}\cap C_{g}\cap \mathbb{Z}^{d}).$
Thus, all cells of the Delaunay tiling defined by $tf+(1-t)g$ are the
intersections of the Delaunay cells of $\mathcal{D}(f)$ and $\mathcal{D}(g);$%
. This implies that $\mathcal{D}(f)$ and $\mathcal{D}(g)$ are commensurate.
\end{proof}

The following Proposition directly follows from Lemma \ref{Commense}.

\begin{proposition}
Two forms $f$ and $g$ belong to the same $L$-cone (not type!) if and only if
their Delaunay tilings $D(f)$ and $D(g)$ are commensurate.
\end{proposition}

\section{G-topes and G*-topes\label{S.G/G*-topes}}

In this section we consider the Delaunay tilings for $\mathbb{Z}^{6}$
relative to the forms $\varphi _{E_{6}}$ and $\varphi _{E_{6}^{\ast }}$.
More specifically, we study the Delaunay tilings of forms lying on the
segment joining $\varphi _{E_{6}}$ and $\varphi _{E_{6}^{\ast }}$of $\frak{P}%
(n)$. These are affinely equivalent to the Delaunay tilings for the root
lattice $E_{6}$ and the dual lattice $E_{6}^{\ast }$. Two Delaunay cells are
said to belong to the same homology class if one of them can be obtained
from the other by the composition of a lattice translation and inversions.

Our discussion is mostly descriptive and adapted specifically for studying
the change of $L$-type along the segment joining forms $\varphi _{E_{6}}$
and $\varphi _{E_{6}^{\ast }}$. More information about these tilings can be
found in Coxeter (1995),\ and Baranovski (1991).

\paragraph{Tiling with G-topes.}

Consider the $45$ $\mathcal{L}_{E_{6}}-$\emph{triangles}$:$\emph{\ } \ \ 
\begin{gather*}
\Delta _{G}^{1}=conv\{[0;0^{5}],[0;1,0^{4}],[2;0,1^{4}]\}\times 5,\quad 
\text{ }\Delta _{G}^{2}=conv\{[3;1^{5}],[-1;-1,0^{4}],[0;1,0^{4}]\}\times 5,
\\
\Delta _{G}^{3}=conv\{[1;0,1^{2},0^{2}],[1;0^{3},1^{2}],[0;1,0^{4}]\}\times
15,\qquad \\
\text{ }\Delta
_{G}^{4}=conv\{[1;1^{2},0^{3}],[-1;0,-1,0^{3}],[2;0,1^{4}]\}\times 20
\end{gather*}
have a common centroid $\mathbf{c}_{G}=\frac{1}{3}[2;1^{5}]$. In these
triples of points admissible positions of non-zero components can be
determined from that the affine centroid is always at $\mathbf{c}_{G}=\frac{1%
}{3}[2;1^{5}]$.\ The convex hull is a polytope with $27$ vertices, the \emph{%
Gosset polytope }$G$ (usually denoted by $2_{21}$), which we call the \emph{%
reference G-tope}. \ There are $54$ homologous copies of $G$ that fit
together facet-to-facet around the origin to form the star at the origin. \
This arrangement can be extended to a tiling by $G$-topes, which is the
Delaunay tiling $\mathcal{D}_{E_{6}}$ for $\mathbb{Z}^{6}$ determined by $%
\varphi _{E_{6}}$ (e.g. Baranovskii (1991)). \ 

These $45$\ triangles inscribed in $G$ are $\mathcal{L}_{E_{6}}-$triangles
because the edge vectors$:$\emph{\ } 
\begin{align*}
E_{G}^{1}& =\{\pm \lbrack 2;-1,1^{4}],\pm \lbrack -2;0,-1^{4}],\pm \lbrack
0;1,0^{4}]\}\times 5, \\
E_{G}^{2}& =\{\pm \lbrack -4;-2,-1^{4}],\pm \lbrack 1;2,0^{4}],\pm \lbrack
3;0,1^{4}]\}\times 5, \\
E_{G}^{3}& =\{\pm \lbrack 0;0,-1^{2},1^{2}],\pm \lbrack
-1;1,0^{2},-1^{2}],\pm \lbrack 1;-1,1^{2},0^{2}]\}\times 15, \\
E_{G}^{4}& =\{\pm \lbrack -2;-1,-2,0^{3}],\pm \lbrack 3;0,2,1^{3}],\pm
\lbrack -1;1,0,-1^{3}]\}\times 20,
\end{align*}
are subsets of $\mathcal{L}_{E_{6}}$. \ This is the set of $270$ long
vectors for $\varphi _{E_{6}}$, and is defined by 
\begin{equation*}
\mathcal{L}_{E_{6}}=\{\mathbf{z\in }\mathbb{Z}^{6}|\varphi _{E_{6}}(\mathbf{z%
})=2m\};
\end{equation*}
$2m$ is the \emph{second minimum} for $\varphi _{E_{6}}$, which is the
minimal value assumed on the non-zero elements of $\mathbb{Z}^{6}$ not
belonging to $\mathcal{S}_{E_{6}}$.\ \ 

As is indicated by the notation, these inscribed\ triangles belong to four $%
\mathcal{A}_{5}-$classes of size $5,5,15$, and $20$, where $\mathcal{A}_{5}$
is the subgroup of permutations of the last five coordinates.\ \ Each long
vector appears only once in the list, so the $45$ edge sets account for the $%
270=45\times 6$ long vectors for $\varphi _{E_{6}}$.

The group $\mathcal{G}_{E_{6}}$\ has the following useful characterization: $%
\mathcal{G}_{E_{6}}\subset $ $GL(6,\mathbb{Z)}$\emph{\ is the stability
group of the set of short vectors }$\mathcal{S}_{E_{6}}$. \ The following
theorem gives information on $\mathcal{G}_{E_{6}}$ (see e.g. Coxeter (1995)
or Baranovskii (1991)).

\begin{theorem}
(Coxeter) \label{T.G-Equivalence}The group $\mathcal{G}_{E_{6}}$ acts
transitively on the $72$ elements of $\mathcal{S}_{E_{6}}$, the $270$
elements of $\mathcal{L}_{E_{6}}$, and, the $54$ G-topes in the star at the
origin.
\end{theorem}

\subsubsection{The inscribed long triangles. \ }

The triangle $\Delta _{G}^{1}$ has a vertex at the origin and is inscribed
in $G$. \ There are five such $\mathcal{L}_{E_{6}}-$triangles that belong to
a single $\mathcal{A}_{5}$ class\ and together have ten vertices at distance 
$2m$ from the origin. \ These vertices, which are given by $\{\mathbf{r}_{1},%
\mathbf{r}_{2},...,\mathbf{r}_{5}\}=\{[0;1,0^{4}]\times 5\},\{\mathbf{b}_{1},%
\mathbf{b}_{2},...,\mathbf{b}_{5}\}=\{[2;0,1^{4}]\times 5\}$, are the
vertices of a five-dimensional cross-polytope. \ The five diagonals
intersect at $\frac{1}{2}[2;1^{5}]$, and the ten diagonal vectors are given
by 
\begin{equation*}
\mathcal{L}_{1}=\{\pm (\mathbf{b}_{1}-\mathbf{r}_{1}),\pm (\mathbf{b}_{2}-%
\mathbf{r}_{2}),...,\pm (\mathbf{b}_{5}-\mathbf{r}_{5})\}=\{\pm \lbrack
2;-1,1^{4}]\times 5\},
\end{equation*}
and belong to a single parity class. \ If $\mathbf{g\in }\mathcal{G}%
_{E_{6}}(G)$, which is the stability group for $G$, then $\mathbf{g}$ must
fix the center $\frac{1}{2}[2;1^{5}]$ and therefore belong to the stability
group of the cross-polytope. \ 

The group of linear transformations that map the cross-polytope onto itself
has order $2^{5}\times 5!$, and can be written as the product $\mathcal{I}%
_{5}\times \mathcal{A}_{5\text{ }}$; the subgroup $\mathcal{I}_{5}$ inverts
arbitrary numbers of axes, and $\mathcal{A}_{5}$ permutes the red vertices $%
\{\mathbf{r}_{1},\mathbf{r}_{2},...,\mathbf{r}_{5}\}$ by permuting the last
five co-ordinates. \ The $\mathcal{A}_{5}-$action clearly stabilizes the set 
$\mathcal{S}_{E_{6}}$, and therefore $\mathcal{A}_{5}\subset \mathcal{G}%
_{E_{6}}(G)$.

The transformations corresponding to

\begin{equation*}
\mathbf{i}_{1}=\frac{1}{2}\left[ 
\begin{array}{rrrrrr}
0 & 4 & 0 & 0 & 0 & 0 \\ 
1 & 0 & 0 & 0 & 0 & 0 \\ 
-1 & 2 & 2 & 0 & 0 & 0 \\ 
-1 & 2 & 0 & 2 & 0 & 0 \\ 
-1 & 2 & 0 & 0 & 2 & 0 \\ 
-1 & 2 & 0 & 0 & 0 & 2
\end{array}
\right] ,\mathbf{i}_{12}=\left[ 
\begin{array}{rrrrrr}
-1 & 2 & 2 & 0 & 0 & 0 \\ 
0 & 0 & 1 & 0 & 0 & 0 \\ 
0 & 1 & 0 & 0 & 0 & 0 \\ 
-1 & 1 & 1 & 1 & 0 & 0 \\ 
-1 & 1 & 1 & 0 & 1 & 0 \\ 
-1 & 1 & 1 & 0 & 0 & 1
\end{array}
\right] ,
\end{equation*}
belong to $\mathcal{I}_{5}$. \ The first inverts the first axis of the
cross-polytope so that 
\begin{equation*}
\mathbf{i}_{1}[\mathbf{r}_{1},\mathbf{r}_{2},...,\mathbf{r}_{5};\mathbf{b}%
_{1},\mathbf{b}_{2},...,\mathbf{b}_{5}]=[\mathbf{b}_{1},\mathbf{r}_{2},...,%
\mathbf{r}_{5};\mathbf{r}_{1},\mathbf{b}_{2},...,\mathbf{b}_{5}],
\end{equation*}
but is non-integral and therefore not an element of $\mathcal{G}_{E_{6}}(G)$%
. \ The second transformation inverts the first two axes so that 
\begin{equation*}
\mathbf{i}_{12}[\mathbf{r}_{1},\mathbf{r}_{2},...,\mathbf{r}_{5};\mathbf{b}%
_{1},\mathbf{b}_{2},...,\mathbf{b}_{5}]=[\mathbf{b}_{1},\mathbf{b}_{2},%
\mathbf{r}_{3},...,\mathbf{r}_{5};\mathbf{r}_{1},\mathbf{r}_{2},\mathbf{b}%
_{3},...,\mathbf{b}_{5}].
\end{equation*}
A direct calculation shows this transformation leaves $\mathcal{S}_{E_{6}}$
invariant, so $\mathbf{i}_{12}\in $ $\mathcal{G}_{E_{6}}(G)$. \ \ Similarly,
there are matrices $\mathbf{i}_{13}$, $\mathbf{i}_{14}$, $\mathbf{i}_{15}\in 
\mathcal{G}_{E_{6}}(G)$ that induce inversions of the pairs of axes
indicated by the subscripts, and together these four inversions generate the
subgroup $\mathcal{I}_{5}^{+}$. \ An arbitrary element $\mathbf{i\in }%
\mathcal{I}_{5}^{+}$ inverts an even number of axes of the cross-polytope. \ 

\begin{proposition}
\label{LongTriangles1.P}$\mathcal{G}_{E_{6}}(G)=\mathcal{I}_{5}^{+}\times 
\mathcal{A}_{5}$, where the $\mathcal{A}_{5}-$action permutes the last four
co-ordinates. $|\mathcal{G}_{E_{6}}(\Delta _{G}^{1})|=2^{4}\times 4!$.
\end{proposition}

\begin{proof}
This follows immediately from the preceeding discussion.
\end{proof}

\begin{corollary}
\label{LongTriangles1.C1}$\mathcal{G}_{E_{6}}$ acts transitively on the $54$
G-topes in the star of G-topes at the origin. \ 
\end{corollary}

\begin{proof}
$|\mathcal{G}_{E_{6}}|/|\mathcal{G}_{E_{6}}(G)|=|\mathcal{G}_{E_{6}}|/|%
\mathcal{I}_{5}^{+}\times \mathcal{A}_{5}|=2^{8}3^{4}5/2^{4}5!=54$. \ \ 
\end{proof}

\noindent

The $\mathcal{G}_{E_{6}}-$action maps the cross-polytope centered at $\frac{1%
}{2}[2;1^{5}]$ to $53$ others, making $54$ it total. \ Since cross-polytopes
are centrally symmetric, opposite $G$-topes have cross-polytopes that are
translates. \ In particular, the cross-polytopes for $G$ and $-G$ are
homologous. \ Consequently, there are $27$ sets of diagonal vectors, each
with $10$ long vectors belonging to the same parity class. \ These sets of
are given by:

\begin{gather*}
\mathcal{L}_{1}=\{\pm \lbrack 2;-1,1^{4}]\times 5\}\times 1;\qquad \mathcal{L%
}_{4}=\{\pm \lbrack 1;2,0^{4}]\times 5\}\times 1; \\
\mathcal{L}_{2}=\{\pm \lbrack -2;0-1^{4}],\pm \lbrack 4;2,1^{4}],\pm \lbrack
0;0,-1^{2},1^{2}],\pm \lbrack 0;0,-1,1,-1,1],\pm \lbrack
0;0,-1,1^{2},-1]\}\times 5; \\
\mathcal{L}_{3}=\{\pm \lbrack 0;1,0^{4}],\pm \lbrack 2;1,2,0^{3}],\pm
\lbrack 2;1,0,2,0^{2}],\pm \lbrack 2;1,0^{2},2,0],\pm \lbrack
2;1,0^{3},2]\}\times 5; \\
\mathcal{L}_{5}=\{\pm \lbrack 3;0,1^{4}],\pm \lbrack -1;0,1,-1^{3}],\pm
\lbrack -1;0,-1,1,-1^{2}], \\
\pm \lbrack -1;0,1^{2},-1,1],\pm \lbrack -1;0,-1^{3},1]\}\times 5; \\
\mathcal{L}_{6}=\{\pm \lbrack -3;-2,0,-1^{3}],\pm \lbrack
-3;0,-2,-1^{3}],\pm \lbrack 1;0^{2},-1,1^{2}], \\
\pm \lbrack 1;0^{2},1,-1,1],\pm \lbrack 1;0^{2},1^{2},-1]\}\times 10.
\end{gather*}
These parity classes are grouped into six $\mathcal{A}_{5}-$classes, as
indicated by the notation. \ This gives a second accounting for the $270$
elements of $\mathcal{L}_{E_{6}}$ in geometric terms.

More directly related to our main line of argument are the following results
on $\Delta _{G}^{1}\subset G$. \ \ 

\begin{corollary}
\label{LongTriangles1.C2}$\ \mathcal{G}_{E_{6}}(\Delta _{G}^{1})=\mathcal{I}%
_{5}^{+}\times \mathcal{A}_{4}$, where $\mathcal{A}_{4}$ is permuation of
the last four coordinates.
\end{corollary}

\begin{proof}
This statement follows from the geometric description of $\Delta _{G}^{1}$
above, and the observation that $\mathcal{G}_{E_{6}}(\Delta _{G}^{1})\subset 
\mathcal{G}_{E_{6}}(G)$.
\end{proof}

\begin{corollary}
\label{LongTriangles1.C3}There are $270$ $\mathcal{L}_{E_{6}}-$triangles $%
\mathcal{G}_{E_{6}}-$equivalent to $\Delta _{G}^{1}$, which are inscribed in
G-topes in the star at the origin.
\end{corollary}

\begin{proof}
It follows from the equality $|\mathcal{G}_{E_{6}}|/|\mathcal{G}%
_{E_{6}}(\Delta _{G}^{1})|=|\mathcal{G}_{E_{6}}|/|\mathcal{I}_{5}^{+}\times 
\mathcal{A}_{4}|=2^{8}3^{4}5/2^{4}4!=270$ that there are $270$ $\mathcal{L}%
_{E_{6}}-$triangles $\mathcal{G}_{E_{6}}-$equivalent to $\Delta _{G}^{1}$. \
That these are inscribed in G-topes in the star follows from \emph{Corollary 
\ref{LongTriangles1.C1}}.
\end{proof}

\paragraph{Tiling with G*-topes.}

The three $\mathcal{L}_{E_{6}^{\ast }}-$\emph{\ triangles } 
\begin{align*}
\Delta _{T}^{1}& =conv\{[0^{6}],[0^{4},1,0],[0^{5},1]\}, \\
\Delta _{T}^{2}& =conv\{[2,0,1^{4}],[-1,0^{2},-1,0^{2}],[-1,0,-1,0^{3}]\}, \\
\Delta _{T}^{3}& =conv\{[0,1,0^{4}],[1,0^{3},1^{2}],[-1^{2},0^{4}]\},
\end{align*}
have a common centroid $\mathbf{c}_{T}=\frac{1}{3}[0^{4},1^{2}]$, and belong
to complementary $2-$spaces. \ Their convex hull is a lattice polytope with
nine vertices, which we denote by $G^{\ast }$. \ This \emph{reference} $%
G^{\ast }-$\emph{tope} is a tile in the Delaunay tiling $\mathcal{D}%
_{E_{6}^{\ast }}$ determined by $\varphi _{E_{6}^{\ast }}$. \ All others are
isometrically equivalent to $G^{\ast }$ relative to $\varphi _{E_{6}^{\ast
}} $. \ In total, there are $40$ homology classes of $G\ast $-topes,
accounting for $720=2\times 9\times 40$ $G\ast $-topes in the star at the
origin. \ 

These inscribed triangles are $\mathcal{L}_{E_{6}^{\ast }}-$ triangles
because the edge vectors 
\begin{align*}
E_{T}^{1}=& \{\pm \lbrack 0^{4},-1,1],\pm \lbrack 0^{5},-1],\pm \lbrack
0^{4},1,0]\}, \\
E_{T}^{2}=& \{\pm \lbrack 0^{2},1,-1,0^{2}],\pm \lbrack -3,0,-2,-1^{3}],\pm
\lbrack 3,0,1,2,1^{2}],\}, \\
E_{T}^{3}=& \{\pm \lbrack 2,1,0^{2},1^{2}],\pm \lbrack -1,-2,0^{4}],\pm
\lbrack -1,1,0^{2},-1^{2}]\},
\end{align*}
belong to $\mathcal{L}_{E_{6}^{\ast }}$, the set of $72$ long vectors for $%
\varphi _{E_{6}^{\ast }}$. \ The long vectors for $\varphi _{E_{6}^{\ast }}$
are defined by 
\begin{equation*}
\mathcal{L}_{E_{6}^{\ast }}=\{\mathbf{z\in }\mathbb{Z}^{6}|\varphi
_{E_{6}^{\ast }}(\mathbf{z})=\frac{3}{2}m\},
\end{equation*}
where $\frac{3}{2}m$ is the second minimum for $\varphi _{E_{6}^{\ast }}$.\ 

There is a second representation $\mathcal{G}_{E_{6}^{\ast }}\subset GL(6,%
\mathbb{Z)}$ for the point group for $E_{6}$, which can be characterized as
either the invariance group for the form $\varphi _{E_{6}^{\ast }}$ or the
stabilizer of $\mathcal{S}_{E_{6}}$. \ The following theorem gives the
essential information we will need on this action (see e.g. Baranovskii
(1991)).

\begin{theorem}
\label{T.G*-Equivalence}The group $\mathcal{G}_{E_{6}^{\ast }}$ acts
transitively on the $54$ elements of $\mathcal{S}_{E_{6}^{\ast }}$, the $72$
elements of $\mathcal{L}_{E_{6}^{\ast }}$, and, the $270$ $G\ast $-topes in
the star at the origin. \ 
\end{theorem}

\noindent From the arithmetic equivalences $\mathcal{P}_{E_{6}^{\ast }}=%
\mathbf{U}_{\ast }\mathcal{S}_{E_{6}^{\ast }}$ and $\mathcal{P}_{E_{6}}=%
\mathbf{U}_{\ast }\mathcal{S}_{E_{6}}$, it follows that $\mathcal{G}%
_{E_{6}}^{\circ }=\mathbf{U}_{\ast }\mathcal{G}_{E_{6}^{\ast }}^{\circ }%
\mathbf{U}_{\ast }^{-1}$ and $\mathcal{G}_{E_{6}^{\ast }}^{\circ }=\mathbf{U}%
_{\ast }\mathcal{G}_{E_{6}}\mathbf{U}_{\ast }^{-1}$, where $\mathbf{U}_{\ast
}\in GL(6,\mathbb{Z})$ was described above. \ The dual groups $\mathcal{G}%
_{E_{6}}^{\circ }=\{\mathbf{g}^{\circ }|\mathbf{g}\in \mathcal{G}_{E_{6}}\}$
and $\mathcal{G}_{E_{6}^{\ast }}^{\circ }=\{\mathbf{g}^{\circ }|\mathbf{g}%
\in \mathcal{G}_{E_{6}^{\ast }}\}$, where $\mathbf{g}^{\circ }=(\mathbf{g}%
^{-1})^{T}$ is the matrix for the dual transformation. \ 

The proof of following corollary will be included into the full-length
version of this paper.

\begin{corollary}
\label{C1.G*-Equivalence}There are $270$ $\mathcal{P}_{E_{6}^{\ast }}-$%
triangles with a vertex at the origin. \ These are $\mathcal{G}%
_{E_{6}}^{\circ }-$equivalent. \ 
\end{corollary}

The $\mathcal{S}_{E_{6}^{\ast }}-$ triangle 
\begin{equation*}
\Delta _{T}^{4}=conv\{[0;0^{5}],[0;1,0^{4}],[2;0,1^{4}]\}
\end{equation*}
is a 2-face of $G^{\ast }$, and is equal to the $\mathcal{L}_{E_{6}}-$
triangle $\Delta _{G}^{1}$. \ The edge set $E_{T}^{4}=$ $E_{G}^{1}=\mathcal{L%
}_{E_{6}}\cap \mathcal{S}_{E_{6}^{\ast }}$. \ 

\begin{corollary}
\label{C2.G*-Equivalence}There are $720$ $\mathcal{S}_{E_{6}^{\ast }}-$
triangles with a vertex at the origin. \ These are $2-$faces of $G^{\ast }$%
-topes, and are $\mathcal{G}_{E_{6}^{\ast }}-$ equivalent.
\end{corollary}

\begin{proof}
This is an immediate consequence of the arithmetic equivalences $\mathcal{S}%
_{E_{6}^{\ast }}=\mathbf{U}_{\ast }^{-1}\mathcal{P}_{E_{6}^{\ast }}$, \ $%
\mathcal{G}_{E_{6}^{\ast }}=\mathbf{U}_{\ast }^{-1}\mathcal{G}_{E_{6}^{\ast
}}\mathbf{U}_{\ast }$, and \emph{Corollary \ref{C1.G*-Equivalence}}.
\end{proof}

\paragraph{Long triangles inscribed in G*-topes.}

The $72$ long vectors for $\varphi _{E_{6}^{\ast }}$ are given by $\mathcal{L%
}_{E_{6}^{\ast }}=\mathcal{L}_{b}\cup \mathcal{L}_{c}$, where $\mathcal{L}%
_{c}=\mathcal{S}_{3}$, and 
\begin{eqnarray*}
\mathcal{L}_{b} &=&\{\pm \lbrack 3,0,1^{4}],\pm \lbrack 1,2,0^{4}],\pm
\lbrack 0^{2};1,0^{3}],\pm \lbrack 2,1;0,1^{3}], \\
&&\pm \lbrack 1,-1;0,1^{3}],\pm \lbrack 3,0;2,1^{3}],\pm \lbrack
1,-1;1^{2},0^{2}]\}
\end{eqnarray*}

The proof of the following lemma will be included into the full-length
exposition of our counterexample.

\begin{lemma}
\label{LongTriangles2.L}There are $120$ homology classes of $\mathcal{L}%
_{E_{6}^{\ast }}-$triangles.
\end{lemma}

\begin{proposition}
\label{LongTriangles2.P}The lattice polytope $T$ is a G*-tope if and only if
it is the convex hull of three $\mathcal{L}_{E_{6}^{\ast }}-$triangles with
a common centroid.
\end{proposition}

\begin{proof}
Each G*-tope is the convex hull of three $\mathcal{L}_{E_{6}^{\ast }}-$%
triangles with a common centroid. \ Since there are $40$ homology classes of
G*-topes, there are $120$ homology classes of $\mathcal{L}_{E_{6}^{\ast }}-$%
triangles inscribed in G*-topes. \ By \emph{Lemma \ref{LongTriangles2.L}}
this accounts for all the homology classes of $\mathcal{L}_{E_{6}^{\ast }}-$%
triangles, from which the proof of the Proposition immediately follows. \ 
\end{proof}

\section{An embedded copy of $D_{4}$\label{S.D4}}

Since $\mathcal{P}_{1}=\{\pm \mathbf{p}_{1},\pm \mathbf{p}_{2},\pm \mathbf{p}%
_{3}\}=\{\pm \lbrack -3;2^{5}],\pm \lbrack 1;0,-1^{4}],\pm \lbrack
2;-2,-1^{4}]\}$ is the edge set for a $\mathcal{P}_{E_{6}^{\ast }}-$%
triangle, this set spans a two-dimensional subspace. \ Let $\Lambda _{D_{4}} 
$ be the four-dimensional sublattice defined by $\mathbb{Z}^{6}\cap \mathcal{%
P}_{1}^{\perp },$and let $\varphi _{D_{4}}$ be the restriction of $\varphi
_{E_{6}}$ to $\Lambda _{D_{4}}$. \ The minimal vectors for $\varphi _{D_{4}}$
are given by $\mathcal{S}_{D_{4}}=\mathcal{S}_{E_{6}}\cap \mathcal{P}%
_{1}^{\perp }=\mathcal{S}_{3}$. \ \ These minimal vectors determine $\varphi
_{D_{4}}$, and $|\mathcal{S}_{D_{4}}|=24$. \ This is sufficient to identify $%
\varphi _{D_{4}}$ is a metrical form for the root lattice $D_{4}$. \ \ 

The set of long vectors for $\varphi _{D_{4}}$ is given by 
\begin{equation*}
\mathcal{L}_{D_{4}}=\mathcal{L}_{E_{6}}\cap \mathcal{P}_{1}^{\perp }=(%
\mathcal{L}_{1}\cap \mathcal{P}_{1}^{\perp })\cup (\mathcal{L}_{2}\cap 
\mathcal{P}_{1}^{\perp })\cup (\mathcal{L}_{3}\cap \mathcal{P}_{1}^{\perp }),
\end{equation*}

which is the union of three parity classes of eight vectors each. \ These
classes are the diagonal vectors for three orientations of cross polytopes,
which tile $\Lambda _{D_{4}}$ by translates. \ This is the Delaunay tiling $%
\mathcal{D}_{D_{4}}$ of $\Lambda _{D_{4}}$ relative to the form $\varphi
_{D_{4}}$. \ This Delaunay tiling can also be realized as the
four-dimensional section $\mathcal{D}_{E_{6}}\cap \mathcal{P}_{1}^{\perp }$.
\ 

\paragraph{The group $\mathcal{G}_{D_{4}}.$}

Let $\mathcal{G}_{D_{4}}=$ $\mathcal{G}_{E_{6}}\cap \mathcal{G}_{E_{6}^{\ast
}}$. \ Since by \emph{Proposition \ref{Invariance.P} }$\mathcal{G}_{E_{6}}$
is the stabilizer of $\mathcal{S}_{E_{6}}=\mathcal{S}_{2}\cup \mathcal{S}$ $%
_{3}$, and $\mathcal{G}_{E_{6}^{\ast }}$ is the stabilizer of $\mathcal{S}%
_{E_{6}^{\ast }}=\mathcal{S}_{1}\cup \mathcal{S}$ $_{2}$, \emph{\ }$\mathcal{%
G}_{D_{4}}$ can be characterized as the subgroup of $GL(6,\mathbb{Z)}$ that
stabilizes each of the sets $\mathcal{S}_{1}$, $\mathcal{S}_{2}$, $\mathcal{S%
}_{3}$. \ 

The short and long vectors for $\varphi _{D_{4}}$ relate to the sets $%
\mathcal{S}_{1}$ and $\mathcal{S}_{3}$ in the following way 
\begin{equation*}
\mathcal{S}_{D_{4}}=\mathcal{S}_{3};\mathcal{L}_{1}\cup \mathcal{L}_{2}\cup 
\mathcal{L}_{3}=\mathcal{S}_{1}\cup \mathcal{L}_{D_{4}},
\end{equation*}
where the last union is disjoint. \ The set $\mathcal{S}_{1}=\{\pm \mathbf{l}%
_{1},\pm \mathbf{l}_{2},\pm \mathbf{l}_{3}\}=\{\pm \lbrack 2;-1,1^{4}],\pm
\lbrack 0;1,0^{4}],$\linebreak $\pm \lbrack -2;0,-1^{4}]\}$ is the edge set $%
E_{G}^{1}$ for the triangle $\Delta _{G}^{1}$, which is inscribed in $G$.

The proof of the following proposition will be included in to a full-length
paper.

\begin{proposition}
\label{D4.P}$\mathcal{G}_{D_{4}}=\mathcal{G}_{E_{6}}(E_{G}^{1})=\mathcal{G}%
_{E_{6}^{\ast }}(E_{G}^{1}).$
\end{proposition}

\medskip

There are six $\mathcal{L}_{E_{6}}-$trianges with a vertex at the origin,
and with edge set $E_{G}^{1}$. \ These fit around the origin, edge-to-edge,
to form a hexagon. \ The triangle $\Delta _{G}^{1}$ is one of the tiles in
this hexagon, and the other five are triangles homologous to it. \ Since $%
\mathcal{G}_{E_{6}}(\Delta _{G}^{1})\subset \mathcal{G}_{E_{6}}(E_{G}^{1})=%
\mathcal{G}_{D_{4}}$, it follows from \emph{Proposition \ref
{LongTriangles1.P}} and \emph{Corollary \ref{LongTriangles1.C1}} that $|%
\mathcal{G}_{D_{4}}|=|\mathcal{G}_{E_{6}}(E_{G}^{1})|=6|\mathcal{G}%
_{E_{6}}(\Delta _{G}^{1})|=2^{8}3^{2}$. \ However, the point group for the
lattice $D_{4}$ has order $2^{7}3^{2}$, which is less than $|\mathcal{G}%
_{D_{4}}|$ by a factor of two. \ This discrepancy is explained by the
following Corollary.

\begin{corollary}
\label{D4.C1}The restriction of $\mathcal{G}_{D_{4}}$ to $\Lambda _{D_{4}}$
is a representation of the point group for the root lattice $D_{4}$. \ 
\end{corollary}

\begin{proof}
Since the point group for $D_{4}$ is the only lattice point group for a
four-dimensional lattice of order $2^{7}3^{2}$, we need only show the kernel
of the restriction map is a two element group. \ The kernal $\mathcal{%
K\subset G}_{D_{4}}$ \ is the subgroup that fixes each point in $\Lambda
_{D_{4}}$.

The four diagonals $(\mathbf{b}_{2}-\mathbf{r}_{2}),...,(\mathbf{b}_{5}-%
\mathbf{r}_{5})$, of the cross-polytope in the long layer of $G$ form a
basis for $\mathcal{P}_{1}^{\perp }$, and therefore an arbitrary element $%
\mathbf{k\in }\mathcal{K}$ must fix these diagonals. \ It follows that $%
\mathbf{k}$ must then fix the parity class $\mathcal{L}_{1}$. These two
conditions can hold only if $\mathbf{k}$ maps $(\mathbf{b}_{1}-\mathbf{r}%
_{1})$ to $\pm (\mathbf{b}_{1}-\mathbf{r}_{1})$. \ By the geometric
discussion of the group $\mathcal{G}_{E_{6}}(\Delta _{G}^{1})$, it follows
that $\mathbf{k\in }$ $\{\mathbf{e},\mathbf{i}_{0}\}\times \mathcal{I}%
_{5}^{+}\subset \mathcal{G}_{D_{4}}$, where $\mathbf{i}_{0}$ is central
inversion and $\mathbf{e}$ is the identity. \ \ Let $\mathbf{k}=\mathbf{k}%
_{1}\mathbf{k}_{2}$, where $\mathbf{k}_{1}\in \{\mathbf{e},\mathbf{i}_{0}\}$%
, $\mathbf{k}_{2}\mathbf{\in }\mathcal{I}_{5}^{+}$. \ If $\ \mathbf{k}_{1}=%
\mathbf{e}$, then $\mathbf{k}_{2}$ must also equal $\mathbf{e}$, and $%
\mathbf{k=e}$. \ If $\mathbf{k}_{1}=\mathbf{i}_{0}$, then $\mathbf{k}_{2}$
must invert the four diagonals $(\mathbf{b}_{2}-\mathbf{r}_{2}),...,(\mathbf{%
b}_{5}-\mathbf{r}_{5})$, which uniquely determines the element $\mathbf{k}%
_{2}$. \ One of the ways to write $\mathbf{k}_{2}$ in this case is as the
product $\mathbf{i}_{23}\mathbf{i}_{34}$, where $\mathbf{i}_{23}$ inverts
the second and third axes of the cross-polytope, and $\mathbf{i}_{45}$
inverts the fourth and fifth. Therefore, when $\mathbf{k}_{1}=\mathbf{i}_{0}$%
, $\mathbf{k=k}_{0}=$ $\mathbf{i}_{0}\mathbf{i}_{23}\mathbf{i}_{34}$, and
the matrix for $\mathbf{k}_{0}$ is given by $\mathcal{S}_{D_{4}}$%
\begin{equation*}
\mathbf{k}_{0}=\left[ 
\begin{array}{rrrrrr}
3 & 0 & -2 & -2 & -2 & -2 \\ 
2 & -1 & -1 & -1 & -1 & -1 \\ 
1 & 0 & 0 & -1 & -1 & -1 \\ 
1 & 0 & -1 & 0 & -1 & -1 \\ 
1 & 0 & -1 & -1 & 0 & -1 \\ 
1 & 0 & -1 & -1 & -1 & 0
\end{array}
\right] .
\end{equation*}
It follows that $\mathcal{K}=\{\mathbf{e},\mathbf{k}_{0}\}$, and is a two
element group.
\end{proof}

\medskip

\paragraph{\ The faces of $\mathcal{D}_{D_{4}}.$}

The origin $\mathbf{0\in \Lambda }_{D_{4}}$, and the star of Delaunay tiles
at the origin has $24$ cross-polytopes, eight belonging to each
orientation.\ \ The lower dimensional cells are $\mathcal{S}_{D_{4}}-$%
simplexes.

\begin{corollary}
\label{D4.C2}For the tiling $\mathcal{D}_{D_{4}}$, the star of $\mathcal{D}%
_{D_{4}}-$cells at the origin contains $24$ edges, $96$ triangles, $96$ $3-$%
simplexes, and $24$ cross-polytopes. \ For each dimension, the cells in the
star are $\mathcal{G}_{D_{4}}-$equivalent. \ 
\end{corollary}

\begin{proof}
By \emph{Corollary \ref{D4.C1}} the homorphic image of $\mathcal{G}_{D_{4}}$
on $\mathcal{\Lambda }_{D_{4}}$ is a group of order $2^{7}3^{2}$, which acts
effectively on $\mathcal{\Lambda }_{D_{4}}$, and is a representation of the
point group for the lattice $D_{4}$. \ The order of the full linear
invariance group of a cross-polytope, attached to the origin in $\mathbb{R}%
^{4}$, is $2^{3}\times 3!$. \ Since $2^{7}3^{2}/24=$ $2^{3}\times 3!$, the
stability group for each of the $24$ cross-polytopes in the star has order $%
2^{3}\times 3!$, and the $\mathcal{G}_{D_{4}}-$orbit of each includes all $%
24 $ cross-polytopes. \ Since the stability group in $\mathcal{G}_{D_{4}}$
of any particular cross-polytope in the star is the full invariance group,
all edges, all $2-$faces and $3-$faces of the cross-polytope, which are
attached to the origin, are $\mathcal{G}_{D_{4}}-$equivalent.

The number of edges attached to the origin is eqaul to the number of short
vectors, which is $24$. \ Each cross-polytope has $12$ $2-$cells attached to
the origin, and each $2-$cell belongs to $3$ cross-polytopes in the star. \
Therefore, there are $(24\times 12)/3=96$ $2-$cells attached to the origin.
\ Each cross-polytope has 8 facets attached to the origin, and each belongs
to a pair of cross-polytopes. \ Therefore, there are $(24\times 8)/2=96$ $3-$%
cells attached to the origin.
\end{proof}

\paragraph{Equivalent facets of $\Phi (\mathcal{P}_{E_{6}^{\ast }})$\ and $%
\Phi (\mathcal{P}_{E_{6}})$.}

\ Since $\mathcal{G}_{E_{6}}^{\circ },\mathcal{G}_{E_{6}^{\ast }}^{\circ
}\subset GL(6,\mathbb{Z)}$ are the full invariance groups of the
corresponding sets of perfect vectors $\mathcal{P}_{E_{6}^{\ast }}$ and $%
\mathcal{P}_{E_{6}}$, $\Phi (\mathcal{P}_{E_{6}^{\ast }})$ is invariant with
respect to the action of $\mathcal{G}_{E_{6}}$, and $\Phi (\mathcal{P}%
_{E_{6}})$ is invariant with respect to the action of $\mathcal{G}%
_{E_{6}^{\ast }}$. \ $\ \ $ The proof of the following corollary will be
included in to a full-length paper.

\begin{corollary}
The perfect domain $\Phi (\mathcal{P}_{E_{6}^{\ast }})$ has $45$ facets that
are $\mathcal{G}_{E_{6}}-$equivalent to $\Phi (\mathcal{P}_{E_{6}^{\ast
}}\cap \mathcal{P}_{E_{6}})$, and, the perfect domain $\Phi (\mathcal{P}%
_{E_{6}})$ has $45$ facets that are $\mathcal{G}_{E_{6}^{\ast }}-$equivalent
to $\Phi (\mathcal{P}_{E_{6}^{\ast }}\cap \mathcal{P}_{E_{6}})$.
\end{corollary}

\noindent The $45$ $\mathcal{G}_{E_{6}}-$equivalent facets of $\Phi (%
\mathcal{P}_{E_{6}^{\ast }})$ described in this Corollary are in one-to-one
correspondence with the $45$ possible edge sets for $\mathcal{L}_{E_{6}}%
\mathcal{-}$triangles inscribed in $G$-topes, and, the $45$ $\mathcal{G}%
_{E_{6}^{\ast }}-$equivalent facets of $\Phi (\mathcal{P}_{E_{6}})$ are in
one-to-one correspondence with the $45$ possible edge sets for $\mathcal{S}%
_{E_{6}^{\ast }}\mathcal{-}$triangles, which are $2-$faces of $G\ast $-topes.

\section{Commensurate and incommensurate G*-topes\label%
{S.IncommensurateTopes}}

A $G^{\ast }$-tope homologous to the reference $G^{\ast }$-tope is called a $%
T$\emph{-tope}. Since $T\subset G$ this $G^{\ast }$-tope is commensurate
with $G$, and therefore commensurate with the Delaunay tiling $\mathcal{D}%
_{E_{6}}$. \ We will count $G^{\ast }$-topes below, and discover that there
are $24$ homology classes of such commensurate $G^{\ast }$-tope. \ 

There is a second \emph{incommensurate} class of $G^{\ast }-$tope with
representative $Q$, which is the convex hull of the three $\mathcal{L}%
_{E_{6}^{\ast }}-$triangles 
\begin{gather*}
\Delta _{Q}^{1}=conv\{[0^{3},1,0^{2}],[0^{4},1,0],[0^{5},1]\},\qquad \Delta
_{Q}^{2}=conv\{[0,1,0^{4}],[-1^{2},0^{4}],[1,0^{2},1^{3}]\}, \\
\Delta _{Q}^{3}=conv\{[2,0,1^{4}],[-1,0^{5}],[-1,0,-1,0^{3}]\},
\end{gather*}
with common centroid $\mathbf{c}_{Q}\mathbf{=}\frac{1}{3}[0^{3},1^{3}]$. \
In the counting exercise below we will find there are $16$ homology classes
of incommensurate $G^{\ast }$-topes, and will refer $G^{\ast }$-topes in any
of these classes as a \emph{Q-topes}. \ There are only two homology classes
of $G^{\ast }$-tope, so the $24$ commensurate classes and the $16$
incommensurate classes sum to the required number of $40$ homology classes
of $G^{\ast }$-tope. the commensurate and incommensurate classes, and the $%
24 $ commensurate classes \ \ There are sixteen homology classes of such $G$%
*-topes, and we will refer to these as \emph{Q-topes}. \ 

An application of \emph{Proposition \ref{LongTriangles2.P}} shows that $Q$
is a $G$*-tope. \ That $Q$ is incommensurate with $G$ is established by
considering the the simplicial $3-$faces $S_{Q}\subset Q$ and $S_{G}\subset
G $. \ The simplex $S_{Q}=conv%
\{[0,1,0^{4}],[2,0,1^{4}],[-1^{2},0^{4}],[-1,0^{5}]\}\subset Q$ is the
convex hull of an edge from each of the triangles $\Delta _{Q}^{2}$, $\Delta
_{Q}^{3}$, so is a $3-$face of $Q$. \ The centroid $\frac{1}{4}[0^{2},1^{4}]$
is also the centroid of the simplex $S_{G}=conv%
\{[0^{2},1,0^{3}],[0^{3},1,0^{2}],[0^{4},1,0],$
$[0^{5},1]\}$, which is a $3-$%
face of $G$. \ These simplexes lie in complementary $3-$spaces, so satisfy
the condition\ that $S_{T}\cap S_{G}=\frac{1}{4}[0^{2};1^{4}]$. \ The
relatively open cells of the intesection polytope $Q\cap G$ are the
intersections of relatively open cells of $Q$ and $G$. \ For this reason the
common centroid is a vertex of $Q\cap G$, and since it is non-integral, $Q$
is incommensurate with $G$. \ Therefore $Q$ is incommensurate with the
Delaunay tiling $\mathcal{D}_{E_{6}}$.\ 

Another important class lattice polytope is the class of \emph{R-topes}; the
lattice polytope $R=conv(S_{Q}\cup S_{G})$ serves as an example. \ This
R-tope has full dimension and eight vertices, so can be triangulated in just
two ways. \ By taking the convex hull of $S_{Q}$ with each $2-$face of $S_{G%
\text{ }}$, the four \emph{blue simplexes} $B^{1},B^{2},B^{3},B^{4}$ are
formed; this collection is the star of simplexes in $R$ with the $3-$face $%
S_{Q}$. \ These blue simplexes tile $R$. The four \emph{yellow simplexes}, $%
Y^{1},Y^{2},Y^{3},Y^{4}$, form a second star, each yellow simplex being the
convex hull of $S_{G}$ with one of the four $2-$faces of $S_{Q}$. \ This is
the collection of simplexes in $R$ with the $3-$face $S_{G}$. \ The yellow
simplexes also tile $R$. \ Lattice polytopes such as $R$ are refered to as 
\emph{repartitioning complexes} in the literature on lattice Delaunay
tilings. \ In the counting exercise below we will see that there are $12$
homology classes of R-topes. \ 

\paragraph{Stars of G-topes.}

In order to effectively study the geometry of Q-topes and R-topes we must
first obtain information on the relationship between the perfect vectors $%
\mathcal{P}_{E_{6}^{\ast }}$, and the G-topes in the star of G-topes at the
origin. \ If $\mathcal{V}(G)$ is the vertex set for $G$, then, relative to $%
\varphi _{E_{6}}$, the set vertices of $G$ that are a distance $m$ from the
origin are given by $\mathcal{S}_{E_{6}}(G)=\mathcal{V}(C)\cap \mathcal{S}%
_{E_{6}}$, and the set vertices that are a distance $2m$ are given by $%
\mathcal{L}_{E_{6}}(C)=\mathcal{V}(C)\cap \mathcal{L}_{E_{6}}$. \ These 
\emph{short} and \emph{long layers} of vertices can also be characterized in
terms of the perfect vector $\mathbf{p}_{1}=[-3,2^{5}]$ using the \emph{%
incidence relations }for G-topes: \ 
\begin{equation*}
\mathcal{S}_{E_{6}}(G)=\{\mathbf{s\in }\mathcal{S}_{E_{6}}|\mathbf{p}%
_{1}\cdot \mathbf{s=}1\};\mathcal{L}_{E_{6}}(G)=\{\mathbf{l\in }\mathcal{L}%
_{E_{6}}|\mathbf{p}_{1}\cdot \mathbf{l=}2\}.
\end{equation*}
The perfect vector $\mathbf{p}_{1}$ determines the \emph{short} and \emph{%
long layers} of vertices, which in turn determine $G$, so we write $G=G_{%
\mathbf{p}_{1}}$. \ This process can be reversed, and the incidence
relations can be used to determine the perfect vector associated with $G$: 
\begin{equation*}
\mathbf{p}_{1}=\mathcal{P}_{G}=\{\mathbf{p\in }\mathcal{P}_{E_{6}^{\ast }}|%
\mathbf{p}\cdot \mathbf{s=}1,\mathbf{s\in }\mathcal{S}_{E_{6}}(G)\}\cap \{%
\mathbf{p\in }\mathcal{P}_{E_{6}^{\ast }}|\mathbf{p}\cdot \mathbf{l=}2,%
\mathbf{l\in }\mathcal{L}_{E_{6}}(G)\}.
\end{equation*}
For this case the set $\mathcal{P}_{G}$ contains the single elements $%
\mathbf{p}_{1}$. \ 

By \emph{Proposition \ref{Invariance.P}}, the incidence relations can be
used to establish a one-to-one correspondence between the $G$-topes in the
star at the origin, and the vectors of $\mathcal{P}_{E_{6}^{\ast }}$: if $%
\mathbf{p\in }\mathcal{P}_{E_{6}^{\ast }}$ then $G_{\mathbf{p}}$ is the
corresponding $G$-tope in the star, and if $G^{\prime }$ is in the star,
then $\mathbf{p}_{G^{\prime }}$ is the corresponding perfect vector. \ This
association can be extended to cover arbitrary cells $C\in \mathcal{D}%
_{E_{6}}$ with a vertex at the origin. \ If $\mathcal{S}_{E_{6}}(C)=\mathcal{%
V}(C)\cap \mathcal{S}_{E_{6}}$ and $\mathcal{L}_{E_{6}}(C)=\mathcal{V}%
(C)\cap \mathcal{L}_{E_{6}}$ are the short and long layers for $C$, define 
\begin{equation*}
\mathcal{P}_{C}=\{\mathbf{p\in }\mathcal{P}_{E_{6}^{\ast }}|\mathbf{p}\cdot 
\mathbf{s=}1,\mathbf{s\in }\mathcal{S}_{E_{6}}(C)\}\cap \{\mathbf{p\in }%
\mathcal{P}_{E_{6}^{\ast }}|\mathbf{p}\cdot \mathbf{l=}2,\mathbf{l\in }%
\mathcal{L}_{E_{6}}(C)\}.
\end{equation*}
When $\dim (C)<6$ the set $\mathcal{P}_{C}$ includes several perfect
vectors. \ The short and long layers, and therefore $C$, can be recoverd
from $\mathcal{P}_{C}$ by again invoking the incidence relations: 
\begin{equation*}
\mathcal{S}_{E_{6}}(C)=\{\mathbf{s\in }\mathcal{S}_{E_{6}}|\mathbf{p}\cdot 
\mathbf{s=}1,\mathbf{p\in }\mathcal{P}_{C}\};\mathcal{L}_{E_{6}}(G)=\{%
\mathbf{l\in }\mathcal{L}_{E_{6}}|\mathbf{p}\cdot \mathbf{l=}2,\mathbf{p\in }%
\mathcal{P}_{C}\}.
\end{equation*}

For a cell $C\in \mathcal{D}_{E_{6}}$, $star_{G}(C)$ is the collection of
G-topes with face $C$. \ The following Proposition gives information on
G-stars when $C$ has a vertex at the origin. \ 

\begin{proposition}
\label{G-star.P}If $\mathbf{0\in }C$, $star_{G}(C)=\{G_{\mathbf{p}}|\mathbf{%
p\in }\mathcal{P}_{C}\}$.
\end{proposition}

\begin{proof}
This assertion follows from the discussion on incidence relations.
\end{proof}

\paragraph{The geometry of R-topes.}

In the statement of the following lemma $star_{G}(S_{G\text{ }})$ is the
collection of G-topes, which have the face $S_{G\text{ }}$, and $%
star_{G^{\ast }}(S_{Q\text{ }})$ is the collection of G*-topes, which have
the face $S_{Q\text{ }}$. \ \ 

The proof of the following lemma will be included into the full-length paper.

\begin{lemma}
\label{R.L}$star_{G^{\ast }}(S_{Q})=\{Q^{1},...,Q^{4}\}$, where $%
Q^{1},...,Q^{4}$ are Q-topes, $Q^{1}=Q$, and the blue simplexes can be
ordered so that $B^{i}=Q^{i}\cap R,i=1,...,4$. \ Similarly, $%
star_{G}(S_{G})=\{G^{1},...,G^{4}\}$, where $G^{1},...,G^{4}$are G-topes, $%
G^{1}=G$, and the yellow simplexes can be ordered so that $Y^{i}=G^{i}\cap
R,i=1,...,4$. \ 
\end{lemma}

\begin{proposition}
\label{R.P}$R$ is commensurate with both $\mathcal{D}_{E_{6}}$ and $\mathcal{%
D}_{E_{6}^{\ast }}$.
\end{proposition}

\begin{proof}
\ \ By \emph{Lemma \ref{R.L}}, the intersection $R\cap \mathcal{D}%
_{E_{6}^{\ast }}$ is tiled by the blue lattice simplexes, so $R$ is
commensurate with $\mathcal{D}_{E_{6}^{\ast }}$. \ Similarly, the
intersection $R\cap \mathcal{D}_{E_{6}}$\ is tiled by the yellow lattice
simplexes, so $R$ is commensurate with $\mathcal{D}_{E_{6}}$.
\end{proof}

\paragraph{The geometry of Q-topes.}

The convex hull of $\Delta _{Q}^{1}$ with one edge from each of the other
triangles $\Delta _{Q}^{2},\Delta _{Q}^{3}\subset Q$ has seven vertices, and
is a lattice simplex in $Q$. \ Two edges can be selected in nine ways, so
there are nine such simplexes, each with $\Delta _{Q}^{1}$ as a $2-$face. \
This is the star of simplexes in $Q$ with the $2-$face $\Delta _{Q}^{1}$,
and this collection tiles $Q$. \ In this subsection we show how six of these
simplexes are commensurate with both $\mathcal{D}_{E_{6}}$ and $\mathcal{D}%
_{E_{6}^{\ast }}$. We call them \emph{white. }The other three are \emph{blue}%
, and can be represented as the intersection of $Q$ with three distinct
R-topes.

\ There are two other collections of simplexes that tile $Q$, the star of
simplexes in $Q$ with the $2-$face $\Delta _{Q}^{2}$, and the star with the $%
2-$face $\Delta _{Q}^{3}$. \ However, both of these tilings play no role in
our discussion. \ The reason for this is that $\Delta _{Q}^{1}$ is a $2-$%
cell in $\mathcal{D}_{E_{6}}$, but the other two triangles $\Delta
_{Q}^{2},\Delta _{Q}^{3}$ are not. \ 

The proof of the following lemma will be included into the full-length paper.

\begin{lemma}
\label{Q.L1}$star_{G}(\Delta _{Q}^{1})=\{G^{1},...,G^{6}\}$, where $%
G^{1},...,G^{6}$ are G-topes and $G^{1}=G$. \ Moreover, for $i=1,...,6$, $%
W^{i}=conv(G^{i}\cap Q\cap \mathbb{Z}^{6})$ is a six-dimensional lattice
simplex with the property that $\Delta _{Q}^{1}\subset W^{i}$. \ We will
refer to these simplexes as white.
\end{lemma}

\begin{lemma}
\label{Q.L2}There are three R-topes $R^{1},R^{2},R^{3}$, such that $R^{1}=R$%
, and such that $B^{i}=R^{i}\cap Q,i=1,...,3$ are six-dimensional lattice
simplexes with the $2-$face $\Delta _{Q}^{1}$. \ Since these simplexes are
the intersections of R-topes and a Q-tope they are blue.
\end{lemma}

\begin{proof}
The three simplexes 
\begin{eqnarray*}
S_{Q}^{1} &=&conv\{[0,1,0^{4}],[2,0,1^{4}],[-1^{2},0^{4}],[-1,0^{5}]\}, \\
S_{Q}^{2} &=&conv\{[0,1,0^{4}],[2,0,1^{4}],[1,0^{2},1^{3}],[-1,0,-1,0^{3}]\},
\\
S_{Q}^{3}
&=&conv\{[-1^{2},0^{4}],[-1,0^{5}],[1,0^{2},1^{3}],[-1,0,-1,0^{3}]\},
\end{eqnarray*}
are $3-$faces of $Q$, and $S_{Q}^{1}=S_{Q}$. \ Each is the convex hull of an
edge from $\Delta _{Q}^{2}$ and $\Delta _{Q}^{3}$. The three simplexes 
\begin{eqnarray*}
S_{G}^{1} &=&conv\{[0^{2},1,0^{3}],[0^{3},1,0^{2}],[0^{4},1,0],[0^{5},1]\},
\\
S_{G}^{2} &=&conv\{[2,1,0,1^{3}],[0^{3},1,0^{2}],[0^{4},1,0],[0^{5},1]\}, \\
S_{G}^{3} &=&conv\{[-2,-1^{2},0^{3}],[0^{3},1,0^{2}],[0^{4},1,0],[0^{5},1]\},
\end{eqnarray*}
are $\mathcal{S}_{D_{4}}-$simplexes, and $S_{G}^{1}=S_{G}$.\ \ The three
centroids of the $Q-$simplexes $S_{Q}^{1},S_{Q}^{2},S_{Q}^{3}$, are given by 
$\frac{1}{4}[0^{2},1^{4}]$, $\frac{1}{4}[2,1,0,2^{3}]$, $\frac{1}{4}[%
-2,-1^{2},1^{3}]$, and these coincide with the centroids of the
corresponding $\mathcal{S}_{D_{4}}-$simplexes $S_{G}^{1},S_{G}^{2},S_{G}^{3}$%
. \ The three $\mathcal{S}_{D_{4}}-$simplexes are translates of closed $3-$%
cells in $\mathcal{D}_{D_{4}}$, so by Proposition \ref{D4.C2} belong to $%
\mathcal{G}_{D_{4}}-$equivalent homology classes. \ It follows that the
three lattice polytopes $R^{1}=conv\{S_{Q}^{1}\cup S_{G}^{1}\}$, $%
R^{2}=conv\{S_{Q}^{2}\cup S_{G}^{2}\}$, $R^{3}=conv\{S_{Q}^{3}\cup
S_{G}^{3}\}$, belong to $\mathcal{G}_{D_{4}}-$equivalent homology classes. \
Since $R^{1}=conv\{S_{Q}^{1}\cup S_{G}^{1}\}=conv\{S_{Q}\cup S_{G}\}=R$,
these three lattice polytopes are R-topes. \ 

Since by Lemma \ref{Q.L2} $B^{1}=Q\cap R^{1}=Q\cap R$ is a six-dimensional
blue lattice simplex, the intersections $B^{2}=Q\cap R^{2}$, $B^{3}=Q\cap
R^{3}$, are other six-dimensional blue lattice simplexes. \ By construction, 
$\Delta _{Q}^{1}\subset B^{1},B^{2},B^{3}$.
\end{proof}

\begin{lemma}
\label{Q.P}The star of simplexes in $Q$ containing $\Delta _{Q}^{1}$, has
six white simplexes $W^{1},W^{2},...,W^{6}$, and three blue simplexes $%
B^{1},B^{2},B^{3}$. \ The six white simplexes are commensurate with $%
\mathcal{D}_{E_{6}}$ and $\mathcal{D}_{E_{6}^{\ast }}$, and the three blue
simplexes are commensurate with $\mathcal{D}_{E_{6}^{\ast }}$, but
incommensurate with $\mathcal{D}_{E_{6}}$.
\end{lemma}

\begin{proof}
The simplexes $W^{1},W^{2},...,W^{6}$, are those described in the statement
of \emph{Lemma \ref{Q.L1}}. \ Since $W^{i}\subset G^{i}\cap Q,i=1,...,6$,
these white simplexes are commensurate with both $\mathcal{D}_{E_{6}}$ and $%
\mathcal{D}_{E_{6}^{\ast }}$. \ 

The simplexes $B^{1},B^{2},B^{3}$, are those described in the statement of 
\emph{Lemma \ref{Q.L2}}. \ Since $B^{i}=R^{i}\cap Q$, the blue simplex $%
B^{i} $ is commensurate with $\mathcal{D}_{E_{6}^{\ast }}$. \ The simplex $%
S_{Q}^{i}$ is a $3-$face of $B^{i}$, and the simplexe $S_{G}^{i}$ is a
closed 3-cells in $\mathcal{D}_{E_{6}}$. \ Since the intersection $%
S_{Q}^{i}\cap S_{G}^{i}$ is the centroid of each of these cells, $S_{Q}^{i}$
and $S_{G}^{i}$, the simplex $B^{i}$ is incommensurate with $\mathcal{D}%
_{E_{6}}$. \ 
\end{proof}

\paragraph{T-topes, Q-topes and R-topes. \ }

The intersection properties of any G*-tope with respect to $\mathcal{D}%
_{E_{6}}$, are invariant with respect to the $\mathcal{G}_{D_{4}}-$action. \
Therefore, any G*-tope $\mathcal{G}_{D_{4}}-$equivalent to $T$ is
commensurate with $\mathcal{D}_{E_{6}}$, and any G*-tope $\mathcal{G}%
_{D_{4}}-$equivalent to $Q$ is incommensurate with $\mathcal{D}_{E_{6}}$. \
Since commensurability is a property that extends to homology classes, it is
natural to make the following definitions: \emph{a T-tope is a lattice
polytope} $\mathcal{G}_{D_{4}}-$\emph{equivalent to one homologous to }$T$%
\emph{\ and, a Q-tope is a lattice polytope }$\mathcal{G}_{D_{4}}-$\emph{%
equivalent to a lattice polytope homologous to }$Q$\emph{\ as a Q-tope.} \
It is also natural to refer to commensurate and incommensurate homology
classes of G*-topes, so that homology classes of T-topes are commensurate
classes and homology classes of Q-topes are incommensurate. \ By the
definitions of T-tope and Q-tope, the homology classes of T-topes are $%
\mathcal{G}_{D_{4}}-$equivalent and, the homology classes of Q-topes are $%
\mathcal{G}_{D_{4}}-$equivalent.

Each of the triangles $\Delta _{T}^{1}$, $\Delta _{T}^{2}$, $\Delta _{T}^{3}$
in $T$ has two edges that are long relative to $\varphi _{E_{6}}$, and a
single edge that is short. \ This is a property that is $\mathcal{G}%
_{D_{4}}- $invariant, so each T-tope has three inscribed triangles of this
type. \ On the other hand, the triangle $\Delta _{Q}^{1}$ in $Q$ has edges
that are short relative to $\varphi _{E_{6}}$, and the triangles $\Delta
_{Q}^{2}$, $\Delta _{Q}^{3}$ have edges that are long. \ Accordingly, each
Q-tope has one short and two long inscribed triangles.

The set of R-topes is the other collection we must consider: \emph{an R-tope
is a lattice polytope} $\mathcal{G}_{D_{4}}-$\emph{equivalent to a lattice
polytope homologous to }$R$. \ The homology classes of R-topes are $\mathcal{%
G}_{D_{4}}-$equivalent, and each R-tope in any of these classes has
identicle intersection properties. \ If $R^{\prime }$ is an arbitrary
R-tope, then $R^{\prime }\cap \mathcal{D}_{E_{6}}$ has four yellow simplexes 
$Y_{R^{\prime \prime }}^{1}$, $Y_{R^{\prime \prime }}^{2}$, $Y_{R^{\prime
\prime }}^{3}$, $Y_{R^{\prime \prime }}^{4}$, which belong to distinct
G-topes and fit around a $3-$cell $S_{G}^{\prime }\in \mathcal{D}_{E_{6}}$; $%
R^{\prime }\cap \mathcal{D}_{E_{6}^{\ast }}$ has four blue simplexes $%
B_{R^{\prime \prime }}^{1}$, $B_{R^{\prime \prime }}^{2}$, $B_{R^{\prime
\prime }}^{3}$, $B_{R^{\prime \prime }}^{4}$, which belong to distinct
Q-topes and fit around a $3-$cell $S_{Q}^{\prime }\in \mathcal{D}%
_{E_{6}^{\ast }}$. \ Each R-tope $R^{\prime }$ is commensurate with both $%
\mathcal{D}_{E_{6}}$ and $\mathcal{D}_{E_{6}^{\ast }}$.

\begin{lemma}
\label{Repackaging.L}The $40$ homology classes of G*-topes are divided
between $24$ homology classes of T-topes, and $16$ homology classes of
Q-topes. \ There are $12$ homology classes of R-topes.
\end{lemma}

\begin{proof}
will be included into the full-length paper
\end{proof}

\paragraph{Repackaging Q-topes.}

Each Q-tope contains six white, and three blue simplexes, and since there
are $16$ homology classes of Q-topes, there are $96$ homology classes of
white, and $48$ homology classes of blue simplexes. \ The $12$ homology
classes of R-topes gives a second accounting for blue simplexes, since each
R-tope can be tiled by four blue simplexes. This allows a repackaging of the
blue simplexes contained in Q-topes, into R-topes. \ The portion of space
tiled by the $16$ classes of Q-topes, can equally well be tiled by the $96$
homology classes of white simplexes, and $12$ classes of R-topes. \ This
alternate tiling has the property that each tile is commensurate with both $%
\mathcal{D}_{E_{6}}$ and $\mathcal{D}_{E_{6}^{\ast }}$.

\begin{proposition}
\label{Structure of D_R}The $24$ homology classes of T-topes, the $12$
homology classes of R-topes, and the $96$ homology classes of white
simplexes fit together facet-to-facet to tile space. This lattice tiling $%
\mathcal{D}_{R}$ is commensurate with both $\mathcal{D}_{E_{6}}$ and $%
\mathcal{D}_{E_{6}^{\ast }}$.
\end{proposition}

\begin{proof}
Since the $T$-topes and $Q$-topes tile space, the discussion immediately
preceeding the proposition shows that the $T$-topes, $R$-topes and white
simplexes also tile space. By \emph{Lemmae} \ref{R.L}, \emph{\ref{Q.P}} and 
\emph{Proposition \ref{R.P}} this tiling is commensurate with both $\mathcal{%
D}_{E_{6}}$ and $\mathcal{D}_{E_{6}^{\ast }}$.
\end{proof}

\section{Proof of main theorem\label{A6.Proof}}

The forms $\varphi _{E_{6}}\in \Phi (\mathcal{P}_{E_{6}^{\ast }})$ and $%
\varphi _{E_{6}^{\ast }}\in \Phi (\mathcal{P}_{E_{6}})$ lie on the central
axes of their domains, and the line segment $\varphi _{t}=(1-t)\varphi
_{E_{6}}+t\varphi _{E_{6}^{\ast }}$, $0\leq t\leq 1$, runs between them. \
The arithmetic minimum along this segment is $m$. \ At the end points the
minimal vectors are respectively $\mathcal{S}_{E_{6}}$ and $\mathcal{S}%
_{E_{6}^{\ast }}$, but at intermediate points the minimal vectors are given
by $\mathcal{S}_{E_{6}}\cap \mathcal{S}_{E_{6}^{\ast }}$.

\paragraph{The form $\protect\varphi _{R}$. \ }

Suppose that the reference $R$ is Delaunay with respect to the metrical form 
$\varphi $. \ Then there is a scalar $c$ and vector $\mathbf{p}$ so that $%
f_{R}(\mathbf{x})=c+\mathbf{p}\cdot \mathbf{x}+\varphi (\mathbf{x})$ is
non-negative on $\mathbb{Z}^{6}$, and zero valued on just the vertices of $R$%
; $f_{R}(\mathbf{x})=0$ is the equation of an empty ellipsoid circumscribing 
$R$. \ The vertex sets for the component simplexes $S_{G}$, $S_{Q}$ of $R$
are given by 
\begin{eqnarray*}
V_{G} &=&\{[0^{2},1,0^{3}],[0^{3},1,0^{2}],[0^{4},1,0],[0^{5},1]\}, \\
V_{Q} &=&\{[0,1,0^{4}],[2,0,1^{4}],[-1^{2},0^{4}],[-1,0^{5}]\},
\end{eqnarray*}
and have the property that $\sum_{\mathbf{v}\in V_{G}}\mathbf{v}=\sum_{%
\mathbf{v}\in V_{Q}}\mathbf{v}=4\mathbf{c}_{R}=[0^{2},1^{4}]$. \ Since $%
f_{R} $ is zero on each vertex of $R$, the metrical form $\varphi $ must
satisfy the condition. 
\begin{eqnarray*}
0 &=&\sum_{\mathbf{v}\in V_{G}}f_{R}(\mathbf{v})-\sum_{\mathbf{v}\in
V_{Q}}f_{R}(\mathbf{v})=\mathbf{p}\cdot \left( \sum_{\mathbf{v}\in V_{G}}%
\mathbf{v-}\sum_{\mathbf{v}\in V_{Q}}\mathbf{v}\right) +\sum_{\mathbf{v}\in
V_{G}}\varphi (\mathbf{v})-\sum_{\mathbf{v}\in V_{Q}}\varphi (\mathbf{v}) \\
&=&\sum_{\mathbf{v}\in V_{G}}\varphi (\mathbf{v})-\sum_{\mathbf{v}\in
V_{Q}}\varphi (\mathbf{v})=trace(\mathbf{P}_{R}\mathbf{F}_{\varphi })=<\pi
_{R},\varphi >,
\end{eqnarray*}
where $\mathbf{F}_{\varphi }$ is the matrix for $\varphi $ and 
\begin{equation*}
\mathbf{P}_{R}=\sum_{\mathbf{v}\in V_{G}}\mathbf{vv}^{T}-\sum_{\mathbf{v}\in
V_{Q}}\mathbf{vv}^{T}=-\left[ 
\begin{array}{rrrrrr}
6 & 1 & 2 & 2 & 2 & 2 \\ 
1 & 2 & 0 & 0 & 0 & 0 \\ 
2 & 0 & 1 & 1 & 1 & 1 \\ 
2 & 0 & 1 & 1 & 1 & 1 \\ 
2 & 0 & 1 & 1 & 1 & 1 \\ 
2 & 0 & 1 & 1 & 1 & 1
\end{array}
\right] .
\end{equation*}
The final expression $<\pi _{R},\varphi >$ is the scalar product of forms $%
\pi _{R}$ and $\varphi $, where $\pi _{R}$ corresponds to $\mathbf{P}_{R}$;
if forms $\pi $, $\varphi $ have coefficient matrices $\mathbf{P}$, $\mathbf{%
F}$, then the scalar product $\left\langle \pi ,\varphi \right\rangle $ is
defined to be $trace(\mathbf{PF)}$. \ 

The equation $<\pi _{R},\varphi >=0$ determines a hyperplane $H_{R}$ in the
space of metrical forms. \ 

\begin{lemma}
\label{Computing the L-type break-point}The line segment $\varphi
_{t}=(1-t)\varphi _{E_{6}}+t\varphi _{E_{6}^{\ast }}$, $0\leq t\leq 1$,
pierces the hyperplane $H_{R}$ at the point $\varphi _{t_{R}}$, where $t_{R}=%
\frac{1}{2}$. \ There are twelve hyperplanes of the form $H_{R^{\prime }}$,
where $R^{\prime }$ is an R-tope, and the line segment $\varphi
_{t}=(1-t)\varphi _{E_{6}}+t\varphi _{E_{6}^{\ast }}$, $0\leq t\leq 1$,
pierces each of these at $\varphi _{1/2}$, which is common to all of them.
\end{lemma}

\begin{proof}
A direct calculation shows that 
\begin{eqnarray*}
\left\langle \pi _{R},\varphi _{t}\right\rangle &=&\left\langle \pi
_{R},(1-t)\varphi _{E_{6}}+t\varphi _{E_{6}^{\ast }}\right\rangle
=(1-t)\times trace(\mathbf{P}_{R}\mathbf{F}_{E_{6}})+t\times trace(\mathbf{P}%
_{R}\mathbf{F}_{E_{6}^{\ast }}) \\
&=&(1-t)(-\frac{1}{2}m)+t(\frac{1}{2}m)=-\frac{1}{2}m+mt.
\end{eqnarray*}
Therefore, the line segment $\varphi _{t},0\leq t\leq 1$, pierces $H_{R}$\
when $t=t_{R}=\frac{1}{2}$.

Each R-tope $R^{\prime }$ corresponds to a form $\pi _{R^{\prime }}$, which
in turn determines a hyperplane $H_{R^{\prime }}$ with equation $<\pi
_{R^{\prime }},\varphi >=0$. \ If $R^{\prime }=$ $\mathbf{g}^{-1}R$, where $%
\mathbf{g\in }\mathcal{G}_{D_{4}}$, then the corresponding form is given by $%
\pi _{R^{\prime }}(\mathbf{x})=\pi _{R}(\mathbf{g}^{\circ }\mathbf{x)}$ and 
\begin{equation*}
H_{R^{\prime }}=\{\varphi (\mathbf{gx)|}\varphi (\mathbf{x)\in }H_{R}\}%
\mathbf{;}
\end{equation*}
$\mathbf{g}^{\circ }$ is the matrix dual to $\mathbf{g}$. \ On the other
hand, if $R^{\prime }$ is homologous to $R$, then $\pi _{R^{\prime }}=\pi
_{R}$ and $H_{R^{\prime }}=H_{R}$. \ \ There are twelve such hyperplanes, as
there are twelve homology classes of R-topes (see \emph{Lemma \ref
{Repackaging.L}}). \ The equality $\varphi _{1/2}(\mathbf{gx)=}\varphi
_{1/2}(\mathbf{x)}$, which holds for all $\mathbf{g\in }\mathcal{G}_{D_{4}},$
implies that the line segment $\varphi _{t}$ pierces all of them at the same
point. \ 
\end{proof}

\begin{proposition}
\label{D_R is Delaunay for Phi}Tiling $\mathcal{D}_{R}$ is Delaunay with
respect to the form $\varphi _{1/2}$.
\end{proposition}

\begin{proof}
Recall that $\mathcal{D}_{R}$ consists of 24 homology classes of $T$-topes,
96 homology classes of white simplexes, and 12 homology classes of $R$-topes
(Proposition \ref{Structure of D_R}). Each $T$\emph{-tope }$T$\emph{\ }in $%
\mathcal{D}_{R}$ is a Delaunay cell of $\mathcal{D}_{E_{6}}$ whose
intersection with a $G$-tope of $\mathcal{D}_{E_{6}^{\ast }}$ is $T.$ Each 
\emph{white} simplex in $\mathcal{D}_{R}$ is the intersection of a $Q$-tope,
Delaunay cell of $\mathcal{D}_{E_{6}^{\ast }},$ and a $G$-tope, Delaunay
cell of $\mathcal{D}_{E_{6}}$. By Lemma \ref{Test on Commensurability of 2
cells} \ all white simplexes and $T$\emph{-topes }are Delaunay with respect
to any intermediate form, including $\varphi _{1/2}.$ By Lemma \ref
{Computing the L-type break-point} and the preceeding discussion, all
homology classes of $R$-topes are Delaunay for $\varphi _{1/2}.$
\end{proof}

\begin{proof}
\label{Commensurate.P}For $0<t<\frac{1}{2}$, the Delaunay tiling for each
form $\varphi _{t}=(1-t)\varphi _{E_{6}}+t\varphi _{E_{6}^{\ast }}$ is
constant, and obtained by replacing each R-tope in $\mathcal{D}_{R}$ by four
yellow simplexes. \ For $\frac{1}{2}<t<1$, the Delaunay tiling for each form 
$\varphi _{t}$ is also constant, and obtained by replacing each R-tope in $%
\mathcal{D}_{R}$ by four blue simplexes.
\end{proof}

\paragraph{The facet $\Phi (\mathcal{P}_{E_{6}}\cap \mathcal{P}_{E_{6}^{\ast
}})$.}

Consider the form $\pi _{P}=\pi _{E_{6}^{\ast }}-\pi _{E_{6}}$ and the rank
one form $\varphi _{\mathbf{p}}(\mathbf{x})=(\mathbf{p\cdot x)}^{2}$. \
Then, $\left\langle \pi _{P},\varphi _{\mathbf{p}}\right\rangle =$ $%
\left\langle \pi _{E_{6}^{\ast }}-\pi _{E_{6}},\varphi _{\mathbf{p}%
}\right\rangle =trace(\mathbf{P}_{E_{6}^{\ast }}-\mathbf{P}_{E_{6}})\mathbf{%
pp}^{T}=\mathbf{p}^{T}(\mathbf{P}_{E_{6}^{\ast }}-\mathbf{P}_{E_{6}})\mathbf{%
p}=\pi _{E_{6}^{\ast }}(\mathbf{p})-\pi _{E_{6}}(\mathbf{p})$. \ If $\mathbf{%
p}\in \mathcal{P}_{E_{6}^{\ast }}$, then $\pi _{P}(\mathbf{p})=\pi
_{E_{6}^{\ast }}(\mathbf{p})-\pi _{E_{6}}(\mathbf{p})=m-\pi _{E_{6}}(\mathbf{%
p})\leq 0$, with equality if an only if $\mathbf{p}\in \mathcal{P}%
_{E_{6}^{\ast }}\cap \mathcal{P}_{E_{6}}$, and, if $\mathbf{p}\in \mathcal{P}%
_{E_{6}}$, then $\pi _{P}(\mathbf{p})=\pi _{E_{6}^{\ast }}(\mathbf{p})-\pi
_{E_{6}}(\mathbf{p})=\pi _{E_{6}^{\ast }}(\mathbf{p})-m\geq 0$, with
equality if an only if $\mathbf{p}\in \mathcal{P}_{E_{6}^{\ast }}\cap 
\mathcal{P}_{E_{6}}$. \ It follows that the hyperplane $H_{P}$, with
equation $\left\langle \pi _{P},\varphi \right\rangle =0$, separates $\Phi (%
\mathcal{P}_{E_{6}})$ and $\Phi (\mathcal{P}_{E_{6}^{\ast }})$, and contains
the facet $\Phi (\mathcal{P}_{E_{6}}\cap \mathcal{P}_{E_{6}^{\ast }})$. \ 

\begin{lemma}
\label{ParameterPerfect} The line segment $\varphi _{t}=(1-t)\varphi
_{E_{6}}+t\varphi _{E_{6}^{\ast }}$, $0\leq t\leq 1$, pierces the facet $%
\mathcal{\Phi (P}_{E_{6}^{\ast }}\cap \mathcal{P}_{E_{6}})$ when $t=t_{P}=%
\frac{2}{5}$. \ Moreover, the point $\varphi _{P}=\varphi _{t_{P}}$ lies on
the central axis of $\mathcal{\Phi (P}_{E_{6}^{\ast }}\cap \mathcal{P}%
_{E_{6}})$.. \ 
\end{lemma}

\begin{proof}
A direct calculation shows that 
\begin{gather*}
\left\langle \pi _{P},\varphi _{t}\right\rangle =\left\langle \pi
_{E_{6}^{\ast }}-\pi _{E_{6}},(1-t)\varphi _{E_{6}}+t\varphi _{E_{6}^{\ast
}}\right\rangle =(1-t)trace(\mathbf{P}_{P}\mathbf{F}_{E_{6}})+(t)trace(%
\mathbf{P}_{P}\mathbf{F}_{E_{6}^{\ast }}) \\
=(1-t)(-\frac{2}{8}m^{2})+t(\frac{3}{8}m^{2})=-\frac{2}{8}m^{2}+\frac{5}{8}%
m^{2}t,
\end{gather*}
where $\mathbf{P}_{P}=\mathbf{P}_{E_{6}^{\ast }}-\mathbf{P}_{E_{6}}$. \
Therefore, the line segment $\varphi _{t},0\leq t\leq 1$, pierces the
hyperplane $H_{P}$ when $t=t_{P}=\frac{2}{5}$. \ 

The forms $\varphi _{t}$ are invariant with respect to the action of $%
\mathcal{G}_{D_{4}}=\mathcal{G}_{E_{6}}\cap \mathcal{G}_{E_{6}}$. \ Since $%
\mathcal{G}_{D_{4}}$\ acts transitively on $\mathcal{P}_{E_{6}}\cap \mathcal{%
P}_{E_{6}^{\ast }}$, which can easily be checked, the only $\mathcal{G}%
_{D_{4}}-$invariant positive forms on $H_{P}$\ must lie on the central axis
of $\Phi (\mathcal{P}_{E_{6}}\cap \mathcal{P}_{E_{6}^{\ast }})$. \
Therefore, the line segment $\varphi _{t}$, $0\leq t\leq 1$, pierces the
hyperplane along the central axis of $\mathcal{\Phi (P}_{E_{6}^{\ast }}\cap 
\mathcal{P}_{E_{6}})$.
\end{proof}

\bigskip

\textbf{Proof of Theorem \ref{Main}. }By Lemma \ref{D_R is Delaunay for Phi} 
$\mathcal{D}_{R}$ is the Delaunay tiling for the midpoint of the segment $%
\varphi _{t}=(1-t)\varphi _{E_{6}}+t\varphi _{E_{6}^{\ast }}.$ By
Proposition \ref{Structure of D_R} $\mathcal{D}_{R}$ is commensurate with
both $\mathcal{D}_{E_{6}}$ and $\mathcal{D}_{E_{6}^{\ast }}$. Therefore by
Lemma \ref{Commense} the $L$-type does not change between $\varphi _{E_{6}}$
and $\varphi _{1/2},$ and between $\varphi _{E_{6}^{\ast }}$ and $\varphi
_{1/2}.$ By the results of Section \ref{FORMS} $\varphi _{E_{6}}$ and $%
\varphi _{E_{6}^{\ast }}$ are the centroids of two \emph{adjacent} perfect
domains and, as shown in the above lemma, the wall between them intersects
the segment $\varphi _{t}$ at $t=\frac{2}{5}.\blacksquare $


\begin{thebibliography}{99}
\bibitem{Alex1}  V. Alexeev, \emph{Complete moduli in the presence of
semiabelian group action }, To appear in \textit{Annals of Mathematics}, 96
pp, arXiv: math.AG/9905103

\bibitem{Alex2}  V. Alexeev, \emph{On Extra Components in the Toroidal
Compactification of $A_{g}$}, To appear in \textit{Moduli of Abelian
Varieties} (Texel island, 1999), Birkhauser, Boston. arXiv: math.AG/9905142

\bibitem{}  E. P. Baranovskii, Partition of Euclidean spaces into $L$%
-polytopes of certain perfect lattices. (Russian) Discrete geometry and
topology (Russian). \textit{Trudy Mat. Inst. Steklov}. \textbf{196} (1991),
27--46; translated in \textit{Proc. Steklov Inst. Math. }\textbf{196, (}%
1992), no. 4, 29-51

\bibitem{}  E. S. Barnes, The complete enumeration of extreme senary forms. 
\textit{Philos. Trans. Roy. Soc. London}. Ser. A. \textbf{249} (1957),
461--506.

\bibitem{B2}  E. S. Barnes, T. J. Dickson, Extreme coverings of $n$-space by
spheres. \textit{J. Austral. Math. Soc}. \textbf{7} (1967), 115--127.

\bibitem{}  Conway, J. H.; Sloane, N. J. A. Low-dimensional lattices. III.
Perfect forms. \textit{Proc. Roy. Soc. London Ser. A} \textbf{418} (1988),
no. 1854, 43--80.

\bibitem{}  Conway, J. H.; Sloane, N. J. A. Low-dimensional lattices. VI.
Voronoi reduction of three-dimensional lattices. \textit{Proc. Roy. Soc.
London Ser. A}\textbf{\ 436} (1992), no. 1896, 55--68.

\bibitem{CS}  Conway, J. H.; Sloane, N. J. A.\textit{\ Sphere packings,
lattices and groups}. Third edition. Grundlehren der Mathematischen
Wissenschaften [Fundamental Principles of Mathematical Sciences], \textbf{290%
}. Springer-Verlag, New York, (1999).

\bibitem{}  H. S. M. Coxeter, \textit{Regular Polytopes}, Dover, NY, 3rd
edition (1973).

\bibitem{Ccl}  Coxeter, H. S. M. Extreme forms.\textit{\ Canadian J. Math}.%
\textbf{\ 3}, (1951). 391--441.

\bibitem{CK}  Coxeter, H. S. M. Kaleidoscopes. \textit{Selected writings of
H. S. M. Coxeter. } CMS Series of Monographs and Advanced Texts. John Wiley $%
\&$ Sons, New York, (1995).

\bibitem{De2}  Delone [Delaunay] B. N. The geometry of positive quadratic
forms, \textit{Uspekhi Mat. Nauk}, \textbf{3} (1937), 16-62; \textbf{4}
(1938), 102-164.

\bibitem{}  Delone, B. N.; Ry\v{s}kov, S. S. Solution of the problem on the
least dense lattice covering of a 4-dimensional space by equal spheres.
(Russian) Dokl. Akad. Nauk SSSR 152, (1963), 523--524.

\bibitem{}  Delone, B. N.; Dolbilin, N. P.; Ry\v{s}kov, S. S.; \v{S}togrin,
M. I. A new construction of the theory of lattice coverings of an \emph{n}%
-dimensional space by congruent balls. (Russian) Izv. Akad. Nauk SSSR Ser.
Mat. 34 1970

\bibitem{DL}  Deza, M. M.; Laurent, M. \textit{Geometry of cuts and metrics.}
Algorithms and Combinatorics, 15. Springer-Verlag, Berlin, 1997.

\bibitem{D2}  Dickson, T. J. A sufficient condition for an extreme covering
of $n$-space by spheres. \textit{J. Austral. Math. So}c. 8 (1968) 56--62.

\bibitem{D1}  Dickson, T. J. On Voronoi reduction of positive definite
quadratic forms.\textit{\ J. Number Theory} 4 (1972), 330--341.

\bibitem{F}  P. Erd\"{o}s, P. M. Gruber, J. Hammer, \textit{Lattice points.}
Pitman Monographs and Surveys in Pure and Applied Mathematics, 39. Longman
Scientific \& Technical, Harlow; copublished in the United States with John
Wiley \& Sons, Inc., New York, 1989.

\bibitem{Holroyd}  Gruber, P. M.; Lekkerkerker, C. G. \textit{Geometry of
numbers}. Second edition. North-Holland Mathematical Library, 37.
North-Holland Publishing Co., Amsterdam-New York, 1987

\bibitem{KZ}  A. Korkine, G. Zolotareff, Sur les formes quadratiques, 
\textit{Math. Ann}. 6 (1873), 366-389.

\bibitem{M1}  J. Martinet, \textit{Les r\'{e}seaux parfaits des espaces
euclidiens.} (French) [Perfect lattices of Euclidean spaces] Math\'{e}%
matiques. [Mathematics] Masson, Paris, 1996

\bibitem{RBcl}  Ryshkov, S. S.; Baranovskii , E. P. $C$-types of $n$%
-dimensional lattices and $5$-dimensional primitive parallelohedra (with
application to the theory of coverings). Cover to cover translation of Trudy
Mat. Inst. Steklov \textbf{137} (1976). Proc. Steklov Inst. Math. 1978, no.%
\textbf{\ 4},\ 140 pp.

\bibitem{}  Ryshkov, S. S.; Baranovskii , E. P., Classical methods in the
theory of lattice packings, \textit{Russian Mathematical Surveys }\textbf{34}%
, 1979, p. 1

\bibitem{}  Voronoi G. F. Nouvelles applications des param\`{e}ters continus 
\`{a} la th\'{e}orie des formes quadratiques, \textit{J. Reine Angew. Math.}%
, Premiere m\'{e}moire, \textbf{133 }(1908)\textbf{, }97-178\textbf{, }Deuxi%
\`{e}me m\'{e}moire\textbf{, 134} (1908), 198-287, \textbf{136} (1909),
67-178.

\bibitem{Ba}  G. F. Voronoi, \textit{Sobranie so\v{c}ineni\u{\i}v treh tomah}%
, [Collected works in three volumes] \textbf{vol. 2}, (in Russian), Kiev
(1952). Introduction and notes by B.N.Delaunay.
\end{thebibliography}
\end{document}